\documentclass[11pt]{article}
\usepackage{amsmath}
\usepackage{amssymb}

\font\tenmsb=msbm10

\def\eps{\varepsilon}
\font\tencmmib=cmmib10 \skewchar\tencmmib '60
\newfam\cmmibfam
\textfont\cmmibfam=\tencmmib

\font\tenmsb=msbm10

\def\Bbb#1{\hbox{\tenmsb#1}}

\def\bbox{\quad\hbox{\vrule \vbox{\hrule \vskip2pt \hbox{\hskip2pt
\vbox{\hsize=1pt}\hskip2pt} \vskip2pt\hrule}\vrule}}
\def\lessim{\ \lower4pt\hbox{$
\buildrel{\displaystyle <}\over\sim$}\ }
\def\gessim{\ \lower4pt\hbox{$\buildrel{\displaystyle >}
\over\sim$}\ }

\def\eps{{\varepsilon}}

\def\Bbb E{\mathbb{E}}
\def\Bbb R{\mathbb{R}} 
\newtheorem{proposition}{Proposition}
\newtheorem{lemma}{Lemma}
\newtheorem{theorem}{Theorem}
\newtheorem{corollary}{Corollary}
\makeatletter
\@addtoreset{equation}{section}

\makeatother

\font\tencmmib=cmmib10 \skewchar\tencmmib '60
\newfam\cmmibfam
\textfont\cmmibfam=\tencmmib

\font\tenmsb=msbm10

\def\Bbb#1{\hbox{\tenmsb#1}}

\def\bbox{\quad\hbox{\vrule \vbox{\hrule \vskip2pt \hbox{\hskip2pt
\vbox{\hsize=1pt}\hskip2pt} \vskip2pt\hrule}\vrule}}
\def\lessim{\ \lower4pt\hbox{$
\buildrel{\displaystyle <}\over\sim$}\ }
\def\gessim{\ \lower4pt\hbox{$\buildrel{\displaystyle >}
\over\sim$}\ }

\def\eps{\varepsilon}

\def\go0{\to 0}

\def\leftitem#1{\item{\hbox to\parindent{\enspace#1\hfill}}}

\def\qed{{$\hfill \bbox$}}

\def\sg{\sigma}

\def\sg2{\sigma^2}

\def\__{_{\infty}}

\begin{document}
{\baselineskip=16.5pt

\title{\bf Von Neumann Entropy Penalization and Low Rank Matrix Estimation}

\author{
{\bf Vladimir Koltchinskii}
\thanks{Partially supported by NSF Grants DMS-0906880 and CCF-0808863}
\\ School of Mathematics
\\ Georgia Institute of Technology
\\ Atlanta, GA 30332-0160
\\ vlad@math.gatech.edu
}

\maketitle

\begin{abstract}
A problem of estimation of a Hermitian nonnegatively 
definite matrix $\rho$ of unit trace 
(for instance, a density matrix of a quantum 
system) based on $n$ independent measurements 
$$
Y_j={\rm tr}(\rho X_j)+\xi_j,\ j=1,\dots, n
$$ 
is studied,
$\{X_j\}$ being i.i.d. Hermitian matrices and $\{\xi_j\}$ 
being i.i.d. mean zero random variables independent of $\{X_j\}.$
 
The estimator
$$
\hat \rho^{\eps}:={\rm argmin}_{S\in {\cal S}}\biggl[
n^{-1}\sum_{j=1}^n (Y_j-{\rm tr}(SX_j))^2 + \eps\ {\rm tr}(S\log S)
\biggr]
$$
is considered, where ${\cal S}$ is the set of all nonnegatively definite Hermitian $m\times m$ matrices of trace $1.$ The goal is to derive oracle inequalities 
showing how the estimation error depends on the accuracy of approximation of the unknown state $\rho$ by low-rank matrices.
\end{abstract}



\medskip

\section{Introduction}

Let ${\mathbb M}_m({\mathbb C})$ be the set of all $m\times m$ matrices 
with complex entries and let  
$$
{\cal S}:=\biggl\{
S\in {\mathbb M}_m({\mathbb C}): S=S^{\ast},\ S \geq 0,
{\rm tr}(S)=1 
\biggr\}
$$ 
be the set of all nonnegatively definite Hermitian matrices of trace $1.$
Here and in what follows $S^{\ast}$ denotes the adjoint matrix of $S$ and 
${\rm tr}(S)$ denotes its trace. 
The matrices from the set ${\cal S}$ can be interpreted, for instance, as \it density matrices, \rm 
describing the states of a quantum system. Given a Hermitian matrix $X$
(\it an observable), \rm its expectation in a state $\rho\in {\cal S}$
is defined as ${\mathbb E}_{\rho}X:={\rm tr}(\rho X).$
Let $X_1,\dots, X_n\in {\mathbb M}_m({\mathbb C}),$ $X_j=X_j^{\ast},\ j=1,\dots, n$ be given Hermitian 
matrices (observables) and let $\rho \in {\cal S}$ be an unknown 
state of the system. An important problem in \it quantum state tomography \rm 
is to estimate $\rho$ based on the observations $(X_j,Y_j),\ j=1,\dots, n,$ where  
$$
Y_j={\rm tr}(\rho X_j)+\xi_j,\ j=1,\dots, n,
$$  
$\xi_j,\ j=1,\dots, n$ being i.i.d. random variables 
with mean zero and finite variance representing measurement errors.
In other words, the unknown state $\rho$ of the system is to be learned 
based on a set of measurements in a number of ``directions''
$X_j, j=1,\dots, n$
(see Artiles, Gill and Guta (2004) for a general discussion 
of statistical problems in quantum state tomography). 
In what follows, it will be usually assumed that the 
design variables $X,X_1,\dots, X_n$ are also random, specifically,
they are i.i.d. Hermitian $m\times m$ matrices with distribution $\Pi,$ and they are independent of the noise $\{\xi_j\}.$  

A typical choice of the design variables 
already discussed in the literature 
(see Gross et al (2009), Gross (2009)) can be described as follows. 
The linear space of matrices ${\mathbb M}_m({\mathbb C})$ can be
equipped with the Hilbert-Schmidt inner product:
$
\langle A,B\rangle := {\rm tr}(AB^{\ast}).
$
Let $E_i,\ i=1,\dots, m^2$ be an orthonormal basis 
of ${\mathbb M}_m({\mathbb C})$ consisting of Hermitian 
matrices $E_i.$ Let $X_j,\ j=1,\dots, n$ 
be i.i.d. random variables sampled from a distribution 
$\Pi$ on the set $\{E_1,\dots, E_{m^2}\}.$
We will refer to this model as 
\it sampling from an orthonormal basis. \rm 
Most often, the uniform distribution $\Pi$ that assigns 
probability $m^{-2}$ to each basis matrix $E_i$ will be used.
Note that in this case 
$
{\mathbb E}|\langle A,X\rangle|^2=m^{-2}\|A\|_2^2, 
$
where $\|\cdot\|_2:=\langle \cdot,\cdot \rangle^{1/2}$ is the Hilbert-Schmidt (or the Frobenius) norm. 

The following simple example is related to the problems of \it matrix completion \rm 
extensively discussed in the recent literature (see, e.g., 
Candes and Recht (2009), Candes and Tao (2009), Recht (2009) and references therein). More precisely, it deals with a version 
of matrix completion for Hermitian matrices (see Gross (2009)).
In this case, when one knows an entry $\rho_{ij}$ of a matrix $\rho,$
one also knows the entry $\rho_{ji}=\bar \rho_{ij}.$

{\bf Example 1. Matrix completion}. Let $\{e_i: i=1,\dots, m\}$ 
be the canonical basis of ${\mathbb C}^m.$ Then, the following set of 
Hermitian matrices forms an orthonormal basis of ${\mathbb M}_m({\mathbb C}):$ 
$$
\Bigl\{e_i \otimes e_i: i=1,\dots, m\Bigr\}
\bigcup 
\biggl\{\frac{1}{\sqrt{2}}(e_i\otimes e_j+e_j\otimes e_i): 1\leq i<j\leq m\biggr\}
$$
$$
\bigcup 
\biggl\{\frac{i}{\sqrt{2}}(e_i\otimes e_j-e_j\otimes e_i): 1\leq i<j\leq m\biggr\},
$$
which will be called \it the matrix completion basis. \rm 
Here and in what follows $\otimes$ denotes the tensor product of 
vectors or matrices. Note that, for a Hermitian matrix $\rho,$ observing inner products $\langle \rho, E_i\rangle$ with randomly picked matrices $E_i$ from the  above basis provides information about real and imaginary parts of the entries of the matrix, which explains the connection to the matrix completion problems. Another option is to consider the following basis 
of the space of all Hermitian matrices:
$$
\Bigl\{e_i \otimes e_i: i=1,\dots, m\Bigr\}
\bigcup 
\biggl\{\frac{1}{2}(e_i\otimes e_j+e_j\otimes e_i)
+\frac{i}{2}(e_i\otimes e_j-e_j\otimes e_i): 1\leq i<j\leq m\biggr\}.
$$
Inner products of a Hermitian matrix $\rho$ with the matrices of this basis are precisely the entries $\rho_{ij}, i\leq j$ of matrix $\rho.$  
If now $\Pi$ is the probability distribution (non-uniform) that assigns probabilities $m^{-2}$ to the matrices $e_i \otimes e_i$ corresponding to the diagonal entries and probabilities $2 m^{-2}$ to other matrices of 
the basis, then 
$
{\mathbb E}|\langle A,X\rangle|^2=m^{-2}\|A\|_2^2.
$ 
Sampling from this distribution is equivalent to sampling the entries 
of the matrix $\rho $ at random (again, recall that when one learns 
an entry $\rho_{ij}$ one also learns $\rho_{ji}=\bar \rho_{ij}$).

Another example was studied by Gross et al (2009) and by Gross (2009). It is more directly related to the problems of quantum state tomography.  

{\bf Example 2. Pauli basis}. Let $m=2^k.$ Consider the \it Pauli basis \rm in the space 
of $2\times 2$ matrices ${\mathbb M}_2({\mathbb C})$: $W_i:=\frac{1}{\sqrt{2}}\sigma_i,$ 
where
$$
\sigma_1:=
\left(
\begin{array}{cc}
0 & 1\\
1 & 0
\end{array}
\right), \ \ \ 
\sigma_2:=
\left(
\begin{array}{cc}
0 & -i\\
i & 0
\end{array}
\right), \ \ \ 
\sigma_3:=
\left(
\begin{array}{cc}
1 & 0\\
0 & -1
\end{array}
\right)\ \ {\rm and}\ \  
\sigma_4:=
\left(
\begin{array}{cc}
1 & 0\\
0 & 1
\end{array}
\right) 
$$
are the \it Pauli matrices. \rm 
Note that the Pauli matrices are both Hermitian and unitary. 
The Pauli basis in ${\mathbb M}_2({\mathbb C})$ can be extended to a basis in the space 
of $m\times m$ matrices ${\mathbb M}_m({\mathbb C}).$ 
These matrices define linear transformations acting in the linear space 
${\mathbb C}^m={\mathbb C}^{2^k}$ that can 
be viewed as a $k$-fold tensor product of spaces 
${\mathbb C}^2:$ ${\mathbb C}^{2^k}=({\mathbb C}^2)^{\otimes k}.$  
Then, the Pauli basis in the space of matrices 
${\mathbb M}_{2^k}({\mathbb C})$ consists of all tensor products $W_{i_1}\otimes \dots \otimes W_{i_k},\ (i_1,\dots, i_k)\in \{1,2,3,4\}^k.$
As before, $X_1,\dots, X_n$ are i.i.d. random variables sampled from 
this set of tensor products. 
Essentially, this is a standard measurement 
model for a $k$ qubit system frequently used in quantum information,
in particular, in quantum state and quantum process tomography
(see Nielsen and Chuang (2000), section 8.4.2).

{\bf Example 3. Subgaussian design}.  
Another interesting class of examples includes \it subgaussian design matrices \rm $X$ such that $\langle A,X\rangle$ is a subgaussian random 
variable for each $A\in {\mathbb M}_m({\mathbb C}).$ (Recall that 
a random variable $\eta$ is called subgaussian with parameter $\sigma$
iff, for all $\lambda \in {\mathbb R},$ ${\mathbb E}e^{\lambda \eta}\leq e^{\lambda^2 \sigma^2/2}$).  
These examples are, probably, of less interest in applications to quantum state tomography, but this is an important 
model, closely related to randomized designs in compressed sensing, 
for which one can use powerful tools developed in the high-dimensional 
probability. For instance, one can consider the \it Gaussian design, \rm where $X$ is a symmetric random matrix with real entries 
such that $\{X_{ij}: 1\leq i\leq j\leq m\}$ are independent centered normal random variables with ${\mathbb E}X_{ii}^2=1,\ i=1,\dots, m$ 
and ${\mathbb E}X_{ij}^2=\frac{1}{2},\ i<j.$   Alternatively, 
one can consider the \it Rademacher design \rm assuming that $X_{ii}=\eps_{ii},\ i=1,\dots, m$ and 
$X_{ij}=\frac{1}{\sqrt{2}}\eps_{ij},\ i<j,$ where $\{\eps_{ij}: 1\leq i\leq j\leq m\}$ are i.i.d. Rademacher random variables (that is, random variables taking values $+1$ or $-1$ with probability $1/2$ each). 
In both cases,
$
{\mathbb E}|\langle A,X\rangle|^2=\|A\|_2^2,\ A\in {\mathbb M}_m({\mathbb C})
$
(such random matrices $X$ will be called \it isotropic\rm)
and $\langle A,X\rangle$ is a subgaussian random variable whose 
subgaussian parameter is equal to $\|A\|_2$ (up to a constant).

The problems of this nature belong to a rapidly growing area of low 
rank matrix recovery. The most popular methods developed so far are based on nuclear norm regularization. 

In what follows, the Euclidean norm in the space ${\mathbb C}^m$ will 
be denoted by $|\cdot |$ and the inner product will be denoted by $\langle \cdot,\cdot \rangle$ (with a little abuse of notation 
since it has been already used for the Hilbert--Schmidt inner 
product between matrices). 
We will denote 
by $\|\cdot\|_{p}, p\geq 1$ the \it Schatten $p$-norm \rm of matrices in ${\mathbb M}_m({\mathbb C})$ (and, if needed, in other matrix spaces). Specifically, 
$
\|A\|_p := \biggl(\sum_{j=1}^m \lambda_k^p (|A|)\biggr)^{1/p},
$
where $|A|:=(A^{\ast}A)^{1/2}$ and, for a Hermitian matrix $B,$
$\lambda_k(B),k=1,\dots, m$ are the eigenvalues of $B$ (usually arranged
in the decreasing order). 
In particular, $\|\cdot\|_1$ is the usual 
nuclear norm and $\|\cdot\|_2$ is the Hilbert-Schmidt norm. 
We will use the notation $\|\cdot\|$ for the operator norm. 
In addition
to the metrics generated by these norms, some other distances will be 
of interest in connection to the statistical problems discussed in 
this paper. In particular, denoting by $\Pi$ the distribution of the design matrix $X,$ we will write 
$$
\|A\|_{L_2(\Pi)}^2:=\int \langle A, x\rangle^2 \Pi(dx)
={\mathbb E}\langle A,X\rangle^2,\ A\in {\mathbb M}_m({\mathbb C})
$$ 
and we will often use the corresponding $L_2(\Pi)$-distance between 
matrices (say, between two states $S_1,S_2\in {\cal S}$). 
This distance represents the prediction error in statistical 
problems in question. 

In the noiseless case (i.e., when $\xi_j \equiv 0$), the following 
estimator of $\rho$ has been extensively studied, especially, in the 
case of matrix completion problems (see Candes and Recht (2009), Candes and Tao (2009), Gross (2009), Recht (2009) and references therein):
$$
\hat \rho :={\rm argmin}\biggl\{\|S\|_1: S\in {\mathbb M}_m({\mathbb C}), \langle S, X_j\rangle=Y_j, j=1,\dots, n\biggr\}.
$$
Under some assumptions that resemble the restricted isometry 
conditions used in compressed sensing, it was shown that, with a high probability, $\hat \rho =\rho$ provided that the number $n$ of observations is sufficiently large. Namely, up to logarithmic factors 
and constants, it should be of the order $m r,$ where $r$ is the 
rank of the target matrix $\rho.$ 

In the noisy case, one has to deal with a matrix regression problem 
and the following penalized least squares estimator,
which is akin to the LASSO used in sparse regression,
was proposed and studied (see, e.g., Candes and Plan (2009), Rohde and 
Tsybakov (2009)):  
\begin{equation}
\label{nuclear_norm}
\hat \rho^{\eps}:={\rm argmin}_{S\in {\mathbb M}_m({\mathbb C})}\biggl[
n^{-1}\sum_{j=1}^n (Y_j-{\rm tr}(SX_j))^2 + \eps \|S\|_1
\biggr],
\end{equation}
where $\eps$ is a regularization parameter. 
Note that these estimators are not 
constrained to the set ${\cal S}$ of density matrices (since for these 
matrices the nuclear norm is equal to $1$).
Candes and Plan (2009)
have also studied another estimator based on the nuclear norm minimization
subject to linear constraints that resembles the Dantzig selector
used in compressed sensing and Rohde and Tsybakov (2009) suggested 
estimators based on nonconvex penalties involving Schatten ``$p$-norms'' 
for $p<1.$ 

We will study the following estimator 
of the unknown state $\rho$ defined as 
a solution of a penalized empirical risk minimization problem: 
\begin{equation}
\label{entropy_penalty}
\hat \rho^{\eps}:={\rm argmin}_{S\in {\cal S}}\biggl[
n^{-1}\sum_{j=1}^n (Y_j-{\rm tr}(SX_j))^2 + \eps\ {\rm tr}(S\log S)
\biggr],
\end{equation}
where $\eps>0$ is a regularization parameter. The penalty term 
is based on the functional 
$
{\rm tr}(S\log S)=-{\cal E}(S),
$
where ${\cal E}(S)$ is the \it von Neumann entropy \rm of state $S.$ 
Thus, the method considered in this paper is based on a trade-off between fitting the model by the least squares in the class of all density matrices and maximizing the entropy of the state. 

One can also consider a slightly different estimator defined as follows:
\begin{equation}
\label{entropy_penalty_2}
\check \rho^{\eps}:={\rm argmin}_{S\in {\cal S}}\biggl[\int \langle S,x\rangle^2 \Pi(dx)
-\frac{2}{n}\sum_{j=1}^n Y_j{\rm tr}(SX_j) + \eps\ {\rm tr}(S\log S)
\biggr].
\end{equation}
Of course, the estimator (\ref{entropy_penalty_2}) requires the knowledge 
of the design distribution $\Pi$ while the estimator (\ref{entropy_penalty}) can be also used in the cases when 
$\Pi$ is unknown. It happens that it is somewhat easier to study the properties of estimator (\ref{entropy_penalty_2}) than of (\ref{entropy_penalty}) for which one has to deal with more complicated 
empirical processes. Note that both optimization problems (\ref{entropy_penalty}) and (\ref{entropy_penalty_2}) are 
convex (this is based on convexity of the penalty term that 
follows from the concavity of von Neumann entropy, see Nielsen and Chuang 
(2000)). In what follows, we will study only the estimators defined 
by (\ref{entropy_penalty}). 

A commutative version of entropy penalization and its connections 
to sparse recovery problems in convex hulls of finite dictionaries 
have been studied by Koltchinskii (2009). In the current paper,  
this approach is extended to the noncommutative case.

\section{An Overview of Main Results}

The results of this paper include oracle inequalities for the $L_2(\Pi)$-error 
of the empirical solution $\hat \rho^{\eps}.$ They will be stated 
in a general form in sections 5 and 6.  
Here we formulate our results only in two of 
the special examples outlined in the Introduction: subgaussian isotropic 
design (such as Gaussian or Rademacher)
and random 
sampling from the Pauli basis.
Assume, for simplicity, that the noise 
$\{\xi_j\}$ is a sequence of i.i.d. $N(0,\sigma_{\xi}^2)$ random variables
(i.e., it is a Gaussian noise).

Let $t>0$ be fixed and denote 
$t_m:=t+\log(2m), \ \ \tau_n:=t+\log\log_2(2n).
$

First we consider the case of subgaussian isotropic design. Note that 
in this case 
$\|A\|_{L_2(\Pi)}=\|A\|_2,\ A\in {\mathbb M}_m({\mathbb C}).$
Given a subspace $L\subset {\mathbb C}^m,$ $P_L$ denotes the 
orthogonal projection on $L$ and $L^{\perp}$ denotes its orthogonal
complement.


\begin{theorem}
\label{intro_th_0}
Suppose $X$ is a subgaussian isotropic matrix.
There exist constants $C>0, c>0$ such that the following 
holds. Under the assumption that $\tau_n\leq c n,$ 
for all $\eps\in [0,1],$
with probability at least $1-e^{-t}$
\begin{eqnarray}
\label{in_A}
&&
\nonumber
\|\hat \rho^{\eps}-\rho\|_{L_2(\Pi)}^2\leq  
C\biggl(
\eps 
\biggl(\|\log \rho\|\wedge 
\log \frac{m}{\eps}\biggr)
\bigvee 
\sigma_{\xi}\sqrt{\frac{m t_m}{n}}\bigvee
\\
&& 
(\sigma_{\xi}\vee \sqrt{m})\frac{\sqrt{m}(\tau_n\log n\vee t_m)}{n} 
\biggr)\biggr].
\end{eqnarray}
Moreover, there exists a constant $D>0$ such that, for all 
$\eps \geq D\sigma_{\xi}\biggl(\sqrt{\frac{m t_m}{n}}\vee \frac{\sqrt{m}t_m}{n}\biggr),$
with probability at least $1-e^{-t},$ 
\begin{eqnarray}
\label{in_B}
&&
\|\hat \rho^{\eps}-\rho\|_{L_2(\Pi)}^2
\leq \inf_{S\in {\cal S}, L\subset {\mathbb C}^m}
\biggl[2\|S-\rho\|_{L_2(\Pi)}^2+
C
\biggl(\eps^2 \|\log S\|_2^2\bigvee 
\sigma_{\xi}^2 \frac{m\ {\rm dim}(L)+\tau_n}{n} \bigvee 
\nonumber
\\
&&
\sigma_{\xi} 
\|P_{L^{\perp}}SP_{L^{\perp}}\|_1
\sqrt{\frac{m t_{m}}{n}}\bigvee
(\sigma_{\xi} \vee \sqrt{m}) \frac{\sqrt{m}(\tau_n \log n\vee t_m)}{n}\biggr)\biggr].
\end{eqnarray}
\end{theorem}

This theorem includes two bounds on the $L_2(\Pi)$-error of 
$\hat \rho^{\eps}.$ The first bound (\ref{in_A}) holds for all $\eps$ including $\eps=0,$ which is the case of the unpenalized least squares estimator. 
The term $\eps \biggl(\|\log \rho\|\wedge \log \frac{m}{\eps}\biggr)$ 
in this bound depends on the operator norm of $\log \rho$ and it has 
to do with the approximation error of the entropy penalization 
method (see Section 4). The second bound (\ref{in_B}) is an oracle inequality 
that controls the squared $L_2(\Pi)$-error of the estimator $\hat 
\rho^{\eps}$ in terms of approximation errors of oracles $S\in {\cal S}.$  
The term $\eps^2 \|\log S\|_2^2$ in this bound is also related
to the approximation error of the entropy penalization method 
discussed in Section 4. This term depends on the 
Hilbert-Schmidt norm of $\log S.$ The dependence on $\eps$ is better 
than in the first bound, but bound (\ref{in_B}) holds only for the values of regularization parameter above certain threshold. 
Clearly, in the second bound, the oracles $S$ are to be of full rank 
(otherwise, $\log S$ does not exist and the right hand side 
of the bound becomes infinite). The random errors in these bounds 
are also different. In the first bound, it is of the order $n^{-1/2}$
(up to logarithmic factors). In the second bound, the error term 
depends on how well the oracle $S$ is approximated by low rank 
matrices. If there exists a subspace $L$ of small dimension ${\rm dim}(L)$ 
such that $\|P_{L^{\perp}}SP_{L^{\perp}}\|_1$ is small (say, 
of the order $n^{-1/2}$), then the random part of the error 
in (\ref{in_B}) is essentially controlled by 
$\sigma_{\xi}^2\frac{{\rm dim}(L) m}{n}.$ 

It will be shown later 
in the paper how to derive from the bounds of 
Theorem \ref{intro_th_0} and more general bounds 
for oracles of full rank some other inequalities for low rank oracles. In particular, for subgaussian isotropic 
design and Gaussian noise, this approach yields the following result. 
To simplify its formulation, we will assume that, for some constant
$c>0,$ $\tau_n\leq cn$ and $t_m\leq n.$


\begin{theorem}
\label{intro_th_2}
Suppose $X$ is a subgaussian isotropic matrix. 
There exist a constant $c>0$ and, for all sufficiently large $D>0,$
a constant $C>0$ such that, for  
$
\eps := D\sigma_{\xi}\sqrt{\frac{m t_m}{n}},
$ 
with probability at least $1-e^{-t},$
\begin{equation}
\label{in_4}
\|\hat \rho^{\eps}-\rho\|_{L_2(\Pi)}^2
\leq \inf_{S\in {\cal S}}\biggl[2 \|S-\rho\|_{L_2(\Pi)}^2+
C\biggl(\frac{\sigma_{\xi}^2{\rm rank}(S) m t_m\log^2(mn)}{n}\bigvee 
\frac{m (\tau_n\log n\vee t_m)}{n}
\biggr)
\biggr].
\end{equation}
\end{theorem}

A simple consequence of the first bound of Theorem \ref{intro_th_0}
and the bound of Theorem \ref{intro_th_2} is 
the following inequality 
that holds with probability at least $1-e^{-t}$ and with some $C>0$
for 
$
\eps := D\sigma_{\xi}\sqrt{\frac{m t_m}{n}}:
$
$$
\|\hat \rho^{\eps}-\rho\|_{L_2(\Pi)}^2 
\leq C\biggl[\biggl(\sigma_{\xi}\sqrt{\frac{m t_m}{n}}\log (mn)
\bigwedge  
\frac{\sigma_{\xi}^2{\rm rank}(\rho) m t_m\log^2(mn)}{n}\biggr)
\bigvee \frac{m (\tau_n \log n \vee t_m)}{n}
\biggr].
$$

\vskip 2mm

Next we consider the case of sampling from the Pauli basis.
In this case, $\|A\|_{L_2(\Pi)}=m^{-1}\|A\|_2,\ A\in {\mathbb M}_m({\mathbb C}).$
As before, we fix $t>0$ and assume that 
$t_m\leq n.$

\begin{theorem}
\label{intro_th_1}
Suppose that $X$ is sampled at random 
from the uniform distribution $\Pi$ on the Pauli basis. 
Then, there exists a constant $C>0$ such that, for all $\eps\in [0,1],$ with probability at least $1-e^{-t},$
\begin{equation}
\label{in_1}
\|\hat \rho^{\eps}-\rho\|_{L_2(\Pi)}^2 
\leq C\biggl[\eps \biggl(\|\log \rho\|\wedge \log\biggl(\frac{m}{\eps}\biggr)\biggr)
\bigvee (\sigma_{\xi}\vee m^{-1/2})\sqrt{\frac{t_m}{nm}}
\biggr].
\end{equation}
In addition, for all sufficiently large $D>0,$ there exists 
a constant $C>0$ such that, for  
$$
\eps := D (\sigma_{\xi}m^{-1/2}\vee m^{-1})\sqrt{\frac{t_m}{n}},
$$
with probability at least $1-e^{-t},$
\begin{equation}
\label{in_2}
\|\hat \rho^{\eps}-\rho\|_{L_2(\Pi)}^2
\leq \inf_{S\in {\cal S}}\biggl[2\|S-\rho\|_{L_2(\Pi)}^2+
C(\sigma_{\xi}^2 \vee m^{-1})\frac{{\rm rank}(S) m t_m\log^2(mn)}{n}\biggr].
\end{equation}
\end{theorem}

Similarly to the previous theorems, one can easily derive 
from Theorem \ref{intro_th_1} the following bound
$$
\|\hat \rho^{\eps}-\rho\|_{L_2(\Pi)}^2 
\leq C\biggl[(\sigma_{\xi}\vee m^{-1/2})\sqrt{\frac{t_m}{mn}}\log (mn)
\bigwedge  
(\sigma_{\xi}^2 \vee m^{-1})\frac{{\rm rank}(\rho) m t_m\log^2(mn)}{n}
\biggr]
$$
that holds with probability at least $1-e^{-t}$ and with some $C>0$
for $\eps = D (\sigma_{\xi}m^{-1/2}\vee m^{-1})\sqrt{\frac{t_m}{n}}.$

It is worth mentioning that the results of sections 4, 5  
provide a way to bound the error of estimator $\hat \rho^{\eps}$
not only in the $L_2(\Pi)$-distance, but also in other statistically 
important distances such as noncommutative Kullback-Leibler, Hellinger and nuclear 
norm distance (see Section 3.1 for their definitions). For instance,
under the assumptions of Theorem \ref{intro_th_0}, the following 
bound for the symmetrized Kullback-Leibler distance
holds with probability at least $1-e^{-t}:$
\begin{eqnarray}
\label{in_Bbb}
&&
K(\hat \rho^{\eps};\rho)
\leq 
\frac{C}{\eps}
\inf_{L\subset {\mathbb C}^m}\biggl[\eps^2 \|\log \rho\|_2^2\bigvee 
\sigma_{\xi}^2 \frac{m\ {\rm dim}(L)+\tau_n}{n} \bigvee 
\nonumber
\\
&&
\sigma_{\xi} 
\|P_{L^{\perp}}\rho P_{L^{\perp}}\|_1
\sqrt{\frac{m t_{m}}{n}}\bigvee
(\sigma_{\xi} \vee \sqrt{m}) \frac{\sqrt{m}(\tau_n \log n\vee t_m)}{n}\biggr].
\end{eqnarray}
In the case of sampling from Pauli basis (as in Theorem \ref{intro_th_1}), it is easy to derive from Theorem \ref{oracle} of Section 5 (using also some bounds from the proofs of Proposition \ref{low_rank_approx} and Corollary \ref{oracle-1})
the following bound on the squared Hellinger distance between 
$\hat \rho^{\eps}$ and $\rho:$
$$
H^2(\hat \rho^{\eps};\rho) 
\leq C
(\sigma_{\xi} \vee m^{-1/2})\frac{{\rm rank}(\rho) 
\sqrt{m t_m}\log^2(mn)}{\sqrt{n}}
$$
that holds with probability at least $1-e^{-t}$ for $\eps = D (\sigma_{\xi}m^{-1/2}\vee m^{-1})\sqrt{\frac{t_m}{n}}.$

It has been already mentioned that the first bounds of theorems \ref{intro_th_0} and \ref{intro_th_1} (bounds (\ref{in_A}) and (\ref{in_1})) hold
for all $\eps\geq 0,$ even in the case of unpenalized least 
squares estimator with $\eps=0.$ The random error parts of these 
bounds are (up to logarithmic factors) of the order $n^{-1/2}$
as $n\to\infty .$ 
Bounds (\ref{in_B}), (\ref{in_4}) and ({\ref{in_2}}) are based on more subtle 
analysis taking into account the ranks of the oracles $S$ approximating 
the true density matrix $\rho .$
In these bounds, the size of the $L_2(\Pi)$-error 
$\|\hat \rho^{\eps}-\rho\|_{L_2(\Pi)}^2$ is determined by a 
trade-off between the approximation error $\|S-\rho\|_{L_2(\Pi)}^2$
of an oracle $S$ and the random error. In the case of bounds (\ref{in_4})
and (\ref{in_2}),  
the last error is of the order 
$\frac{\sigma_{\xi}^2{\rm rank}(S)m}{n}$ 
(up to logarithmic factors), 
and 
it depends on the rank of the oracle $S.$ In particular, taking 
$S=\rho,$ we can conclude that $\|\hat \rho^{\eps}-\rho\|_{L_2(\Pi)}^2$
is bounded by $\frac{\sigma_{\xi}^2{\rm rank}(\rho)m}{n}$ (up to constants and logarithmic factors). This means that von Neumann entropy penalization mimics oracles that know precisely which low rank matrices approximate 
$\rho$ well and can estimate $\rho$ by estimating a ``small'' number 
of parameters needed to describe such oracles. This could be compared 
with recent results for nuclear norm penalization (Candes and Plan (2009), Rohde and Tsybakov (2009)). 
Depending on the values 
of $\sigma_{\xi}, m, n$ and other characteristics of the problem more 
``rough'' bounds 
(\ref{in_A}) and ({\ref{in_1}}) might become even sharper 
than more ``subtle'' bounds (\ref{in_B}), ({\ref{in_4}}) and ({\ref{in_2}})
(see Rohde and Tsybakov (2009) for a discussion of a similar 
phenomenon). Since the random error term in more ``subtle'' 
bounds is proportional to $\sigma_{\xi}^2$ and in the ``rough'' 
bounds it is proportional to $\sigma_{\xi},$ the ``rough'' 
bounds become sharper for the values of standard deviation of the noise $\sigma_{\xi}$ above a threshold that depends on $n$ and $m.$  
Thus, the rate of convergence of the $L_2(\Pi)$-error 
to zero in a particular asymptotic scenario (when certain characteristics 
are large) is determined by the bounds of both types.

Theorems \ref{intro_th_0}, \ref{intro_th_2}, \ref{intro_th_1} and other results of a similar nature will follow as corollaries from more general oracle inequalities that we establish under broader assumptions on the design distributions and on the noise. 
To prove these results, we need several tools from  the empirical processes and random matrices theory, such as noncommutative Bernstein type inequalities and 
generic chaining bounds for empirical processes.
We will discuss these results in Section 3 (as well as some properties 
of noncommutative Kullback-Leibler, Hellinger and other distances between 
density matrices). We will then study approximation error bounds for 
the solution of von Neumann entropy penalized true risk minimization 
problem (Section 4) and, finally, in sections 5 and 6, derive main 
results of the paper concerning random error bounds for the empirical
solution $\hat \rho^{\eps}.$ More precisely, we bound the squared $L_2(\Pi)$-distance $\|\hat \rho^{\eps}-S\|_{L_2(\Pi)}^2$ and symmetrized Kullback-Leibler distance $K(\hat \rho^{\eps};S)$ from $\hat \rho^{\eps}$
to an arbitrary ``oracle'' $S\in {\cal S}$  and derive oracle inequalities for the squared $L_2(\Pi)$-error $\|\hat \rho^{\eps}-\rho\|_{L_2(\Pi)}^2$ of the empirical solution $\hat \rho^{\eps}.$ These results are first established for oracles $S$ of full rank and expressed in terms of certain
characteristics of the operator $\log S$ (which is, essentially, a subgradient of the von Neumann entropy penalty used in (\ref{entropy_penalty})). 
Using simple techniques discussed in Section 4, we then develop 
the bounds for low rank oracles $S$ (such as the bounds of theorems \ref{intro_th_2} and \ref{intro_th_1}) and also obtain oracle inequalities 
for so called ``Gibbs oracles''. Note that the logarithmic factors involved in the bounds of theorems \ref{intro_th_2} and \ref{intro_th_1}
(and in other results of this type discussed later in the paper),
in particular, the factor $\log^2(mn),$ are related to the need 
to bound certain norms of $\log S$ for special oracles $S\in {\cal S}$
(as in Theorem \ref{intro_th_0}). 
In the case of $\|S\|_1$-penalization, $\log S$ should be replaced with 
a version of ${\rm sign}(S)$ and one can avoid 
some of the logarithmic factors in this case.

\section{Preliminaries: Distances in ${\cal S},$ Empirical Processes and Exponential Inequalities for Random Matrices}

\subsection{Noncommutative Kullback-Leibler and other distances} 

We will use noncommutative extensions of 
classical distances between probability distributions such as 
Kullback-Leibler and Hellinger distances. These extensions 
are common in quantum information theory (see Nielsen and Chuang (2000)). 
In particular, we will use \it Kullback-Leibler divergence \rm between 
two states $S_1,S_2\in {\cal S}$ defined as 
$$
K(S_1\|S_2):={\mathbb E}_{S_1}(\log S_1-\log S_2)={\rm tr}(S_1(\log S_1-\log S_2))
$$
and its symmetrized version 
$$
K(S_1;S_2):=K(S_1\|S_2)+K(S_2\|S_1)={\rm tr}((S_1-S_2)(\log S_1-\log S_2)).
$$
We will also use a noncommutative version of \it Hellinger distance \rm defined 
as follows. For any two states $S_1,S_2\in {\cal S},$ let
$
F(S_1,S_2):={\rm tr}\sqrt{S_1^{1/2}S_2S_1^{1/2}}.
$ 
This quantity is called the \it fidelity \rm of states $S_1,S_2$
(see, e.g., Nielsen and Chuang (2000), p. 409). Then, a natural 
definition of the squared Hellinger distance is 
$
H^2(S_1,S_2):=2(1-F(S_1,S_2)).
$ 
A remarkable property of this distance is that 
$$
H^2(S_1,S_2)=\sup H^2(\{p_i\};\{q_i\})=\sup \sum_{i} \Bigl(\sqrt{p_i}-\sqrt{q_i}\Bigr)^2,
$$
where the supremum is taken over all POVMs $\{E_i\}$ 
(positive operator valued measures) and $p_i:={\rm tr}(S_1 E_i),
q_i:={\rm tr}(S_2 E_i).$ [In the discrete case, a positive operator 
valued measure is a set $\{E_i\}$ of Hermitian nonnegatively definite 
matrices such that $\sum_{i}E_i=I$]. 
Thus, the quantum Hellinger distance 
is just the largest ``classical'' Hellinger distance between
the probability distributions $\{p_i\}, \{q_i\}$ of a ``measurement''
$\{E_i\}$ in the states $S_1,S_2$ (see Nielsen and Chuang (2000), p. 412).
The same property also holds for two other important ``distances'',
the trace distance $\|S_1-S_2\|_1$ and the Kullback-Leibler 
divergence $K(S_1\|S_2)$ (see, e.g., Klauck et al (2007)).
These properties immediately imply
an extension of classical inequalities for these distances: 
$$
\|S_1-S_2\|_1^2 \leq H^{2}(S_1,S_2)\leq K(S_1\|S_2).
$$
They also imply the following simple proposition used below.
It shows that, if two matrices $S_1,S_2$ are close in the Hellinger 
distance and one of them (say, $S_2$) is ``approximately low rank'' 
in the sense that there exists a subspace $L\subset {\mathbb C}^m$ of small dimension 
such that $\|P_{L^{\perp}}S_2P_{L^{\perp}}\|_1$ is small,
then another matrix $S_1$ is also ``approximately low rank''
with the same ``support'' $L.$ 

\begin{proposition}
\label{rank_00}
For all subspaces $L\subset {\mathbb C}^m$ and all $S_1,S_2\in {\cal S},$
$$
\|P_L S_1 P_L\|_1 
\leq 2\|P_L S_2 P_L\|_1 + 2 H^2(S_1,S_2).
$$ 
\end{proposition}
 
{\bf Proof}. Indeed, take an orthonormal basis $\{e_1,\dots, e_m\}$ in ${\mathbb C}^m$
such that $L={\rm l.s.}(\{e_1,\dots, e_k\}).$ Let 
$p_j:=\langle S_1 e_j,e_j\rangle= {\rm tr}(S_1 (e_j\otimes e_j))$ 
and   
$q_j:=\langle S_2 e_j,e_j\rangle= {\rm tr}(S_2(e_j\otimes e_j)).$ 
Then 
$$
H^{2}(S_1,S_2)\geq 
\sum_{j=1}^m \Bigl(\sqrt{p_j}-\sqrt{q_j}\Bigr)^2\geq 
\sum_{j=1}^k \Bigl(\sqrt{p_j}-\sqrt{q_j}\Bigr)^2
=\sum_{j=1}^k p_j + 
\sum_{j=1}^k q_j - 2\sum_{j=1}^k \sqrt{p_j}\sqrt{q_j},
$$
which implies (using that $2\sqrt{ab}\leq a/2+2b$)
$$
\|P_L S_1 P_L\|_1 =
\sum_{j=1}^k p_i
\leq 
2\sum_{j=1}^k \sqrt{p_j}\sqrt{q_j}
-\sum_{j=1}^k q_j + H^{2}(S_1,S_2)
\leq 
$$
$$
\frac{1}{2}\sum_{j=1}^k p_j+ \sum_{j=1}^k q_j 
+H^2(S_1,S_2)=
\frac{1}{2}\|P_L S_1 P_L\|_1 
+\|P_L S_2 P_L\|_1 + H^2(S_1,S_2),
$$
and the result follows.

\qed

\subsection{Empirical processes bounds}

We will use several inequalities for empirical processes indexed 
by a class of measurable functions ${\cal F}$ defined on an arbitrary 
measurable space $(S,{\cal A}).$ Let $X, X_1,\dots ,X_n$ be i.i.d. random 
variables in $(S,{\cal A})$ with common distribution $P.$ If ${\cal F}$
is uniformly bounded by a number $U,$ then Bousquet's version of the 
famous Talagrand's concentration inequality for empirical processes implies that, for all $t>0,$ with probability at least $1-e^{-t}$
$$
\sup_{f\in {\cal F}}\biggl|n^{-1}\sum_{j=1}^n f(X_j)-{\mathbb E}f(X)\biggr|
\leq 2\biggl[{\mathbb E}\sup_{f\in {\cal F}}\biggl|n^{-1}\sum_{j=1}^n f(X_j)-{\mathbb E}f(X)\biggr|+
\sigma \sqrt{\frac{t}{n}}+U\frac{t}{n}\biggr],
$$
where 
$
\sigma^2 := \sup_{f\in {\cal F}}{\rm Var}_P(f(X)).
$
We will also need a version of this bound for function classes that are 
not necessarily uniformly bounded. Such a bound was recently proved by 
Adamczak (2008). Let 
$
F(x)\geq \sup_{f\in {\cal F}}|f(x)|, x\in S,
$
be an envelope of the class. It follows from Theorem 4 of Adamczak (2008)
that, there exists a constant $K>0$ such that for all $t>0$ with probability
at least $1-e^{-t}$$$
\sup_{f\in {\cal F}}\biggl|n^{-1}\sum_{j=1}^n f(X_j)-{\mathbb E}f(X)\biggr|
\leq K\biggl[
{\mathbb E}\sup_{f\in {\cal F}}\biggl|n^{-1}\sum_{j=1}^n f(X_j)-{\mathbb E}f(X)\biggr|+
\sigma \sqrt{\frac{t}{n}}+\Bigl\|\max_{1\leq j\leq n}|F(X_j)|\Bigr\|_{\psi_{1}}\frac{t}{n}\biggr].
$$

In addition to this, we will need to bound the following expectation:
$$
{\mathbb E}\sup_{f\in {\cal F}}\biggl|n^{-1}\sum_{j=1}^n f^2(X_j)-{\mathbb E}f^2(X)\biggr|.
$$
A usual approach to this problem is to use symmetrization inequality 
to replace the empirical process by a Rademacher process, and then 
to use Talagrand's comparison (contraction) inequality (see, e.g., 
Ledoux and Talagrand (1991), Section 4.5) to get rid of the squares. This, however,
would require the class ${\cal F}$ to be uniformly bounded by some $U>0,$
which is not too large. This approach is not sufficient in the case 
of subgaussian design considered in the last section. A more subtle 
approach has been developed in the recent years by Klartag and 
Mendelson (2005), Mendelson (2010) and it is based on generic 
chaining bounds.

Talagrand's \it generic chaining complexity \rm 
(see Talagrand (2005)) of a metric space 
$(T,d)$ is defined as follows. An admissible sequence $\{\Delta_n\}_{n\geq 0}$ is an increasing sequence of partitions of $T$ (i.e., each next 
partition is a refinement of the previous one) such that ${\rm card}(\Delta_0)=1$ and ${\rm card}(\Delta_n)\leq 2^{2^n},\ n\geq 1.$
For $t\in T,$ $\Delta_n(t)$ denotes the unique subset in $\Delta_n$
that contains $t.$ For a set $A\subset T,$ $D(A)$ denotes its diameter. 
Then, define the generic chaining complexity $\gamma_2(T;d)$ as 
$$
\gamma_2(T;d):=\inf_{\{\Delta_n\}_{n\geq 0}} \sup_{t\in T}\sum_{n\geq 0} 2^{n/2} D(\Delta_n(t)),
$$
where the $\inf $ is taken over all admissible sequences of partitions.  

If $\{X(t): t\in T\}$ is a centered Gaussian process such that 
$
{\mathbb E}(X(t)-X(s))^2 = d^2 (t,s),\ t,s\in T, 
$
then it was proved by Talagrand that 
$$
K^{-1}\gamma_2(T;d)\leq {\mathbb E}\sup_{t\in T}X(t)\leq K \gamma_2(T;d),
$$
where $K>0$ is a universal constant. Thus, the generic chaining complexity
$\gamma_2(T;d)$ is a natural characteristic of the size of the Gaussian 
process $X(t), t\in T.$ 

Similar quantities can be also used to control 
the size of empirical processes indexed by a function class ${\cal F}.$
It is natural to define $\gamma_2({\cal F};L_2(P)),$ that is, $\gamma_2({\cal F};d),$ where $d$ is the $L_2(P)$-distance. Some other distances are also useful, for instance, the $\psi_2$-distance associated 
with the probability space $(S,{\cal A}, P).$
Recall that, for a convex increasing function $\psi$ with $\psi(0)=0,$ 
$$
\|f\|_{\psi}:=\inf\biggl\{C>0: \int_S\psi\biggl(\frac{|f|}{C}\biggr)dP\leq 1\biggr\}
$$
(see van der Vaart and Wellner (1996), p. 95).
If $\psi (u)=u^p, u\geq 0,$ for some $p\geq 1,$ the corresponding 
$\psi$-norm is just the $L_p$-norm. Other important choices
are functions $\psi_{\alpha}(t)=e^{t^{\alpha}}-1, t\geq 0, \alpha\geq 1,$
especially, $\psi_2$ that is related to subgaussian tails of $f$ and 
$\psi_1$ that is related to subexponential tails.

The generic chaining complexity that corresponds to the $\psi_2$-distance  will be denoted 
by $\gamma_2({\cal F};\psi_2).$ Mendelson (2010) proved the following
deep result (strengthening previous results by Klartag and Mendelson (2005)). Suppose that ${\cal F}$ is a symmetric class, that is, 
$f\in {\cal F}$ implies $-f\in {\cal F},$ and $Pf={\mathbb E}f(X)=0, f\in {\cal F}.$ Then, for some universal constant $K>0,$
$$
{\mathbb E}\sup_{f\in {\cal F}}\biggl|n^{-1}\sum_{j=1}^n f^2(X_j)-{\mathbb E}f^2(X)\biggr|
\leq K\biggl[
\sup_{f\in {\cal F}}\|f\|_{\psi_1}\frac{\gamma_2({\cal F};\psi_2)}{\sqrt{n}}
\bigvee \frac{\gamma_2^2({\cal F};\psi_2)}{n}
\biggr].
$$

\subsection{Noncommutative Bernstein type inequalities}

We will need the following operator version of Bernstein's 
inequality which is due to Ahlswede and Winter (2002) 
(and which has been already successfully used in the low rank 
recovery problems by Gross et al (2009), Gross (2009), Recht (2009)). 

In this subsection, assume that 
$X,X_1,\dots, X_n$ are i.i.d. random Hermitian $m\times m$
matrices with ${\mathbb E}X=0$ and $\sigma_X^2:=\|{\mathbb E}X^2\|.$

{\bf Bernstein's inequality for operator valued r.v.}
{\it Suppose that 
$\|X\|\leq U$ for some $U>0.$ 
Then 
\begin{equation}
\label{bernstein-0}
{\mathbb P}\biggl\{\|X_1+\dots+X_n\|\geq t\biggr\}
\leq 2m \exp\biggl\{-\frac{t^2}{2\sigma_X^2 n+2Ut/3}\biggr\}.
\end{equation}
}

In fact, we will frequently use the following bound that immediately 
follows from the version of Bernstein's inequality given above:
{\it  
for all $t>0,$ with probability at least $1-e^{-t}$
\begin{equation}
\label{bernstein-1}
\biggl\|\frac{X_1+\dots+X_n}{n}\biggr\|\leq 
2\biggl(\sigma_X\sqrt{\frac{t+\log(2m)}{n}}\bigvee U\frac{t+\log(2m)}{n}\biggr). 
\end{equation}
}

Moreover, it is possible to replace the $L_{\infty}$-bound $U$ on 
$\|X\|$ in the above inequality by bounds on the weaker $\psi_{\alpha}$-norms. 
Denote 
$
U_X^{(\alpha)}:=\Bigl\|\|X\|\Bigr\|_{\psi_{\alpha}},\ \alpha\geq 1.
$

\begin{proposition} 
\label{bernstein_A}
Let $\alpha\geq 1.$ There exists a constant $C>0$ such that, 
for all $t>0,$ with probability at least $1-e^{-t}$
\begin{equation}
\label{bernstein-2}
\biggl\|\frac{X_1+\dots+X_n}{n}\biggr\|\leq 
C\biggl(\sigma_X\sqrt{\frac{t+\log(2m)}{n}}\bigvee 
U_{X}^{(\alpha)}\biggl(\log \frac{U_X^{(\alpha)}}{\sigma_X}\biggr)^{1/\alpha}\frac{t+\log(2m)}{n}\biggr). 
\end{equation}
\end{proposition}
 
Note that, in the limit $\alpha \to \infty,$ inequality (\ref{bernstein-2}) coincides with (\ref{bernstein-1}) (up to 
a constant). 

{\bf Proof}. 
Similarly to the proof of (\ref{bernstein-0}) discussed in the 
literature (Ahlswede and Winter (2002), Gross (2009), Recht (2009)),
we follow the standard derivation of classical Bernstein's inequality
and we use the well known \it Golden-Thompson inequality \rm (see, e.g.,  
Simon (1979), p. 94): for arbitrary Hermitian matrices $A,B\in {\mathbb M}_m({\mathbb C}),$
$
{\rm tr}(e^{A+B})\leq {\rm tr}(e^A e^B).
$ 
Let $Y_n:=X_1+\dots +X_n.$ Note that $\|Y_n\|< t$ if and only if 
$-tI_m<Y_n<tI_m.$ Therefore,
\begin{equation}
\label{odin_1}
{\mathbb P}\{\|Y_n\|\geq t\}={\mathbb P}\{Y_n\not\leq tI_m\}+
{\mathbb P}\{Y_n\not\geq -tI_m\}.
\end{equation}
The following bounds are straightforward by simple matrix algebra:
\begin{equation}
\label{dwa_2}
{\mathbb P}\{Y_n \not\leq tI_m\}=
{\mathbb P}\{e^{\lambda Y_n}\not\leq 
e^{\lambda t I_m}\}\leq 
{\mathbb P}\Bigl\{{\rm tr}\Bigl(e^{\lambda Y_n}\Bigr)\geq 
e^{\lambda t }\Bigr\}
\leq e^{-\lambda t}
{\mathbb E}{\rm tr}(e^{\lambda Y_n}).
\end{equation}
To bound the expected value in the right hand side, we use 
independence of random variables $X_1,\dots, X_n$ and Golden-Thompson inequality:
$$
{\mathbb E}{\rm tr}(e^{\lambda Y_n})=
{\mathbb E}{\rm tr}\Bigl(e^{\lambda Y_{n-1}+\lambda X_n}\Bigr)
\leq 
{\mathbb E}{\rm tr}\Bigl(e^{\lambda Y_{n-1}}e^{\lambda X_n}\Bigr)=
{\rm tr}\biggl({\mathbb E}\Bigl(e^{\lambda Y_{n-1}}e^{\lambda X_n}\Bigr)\biggr)=
$$
$$
{\rm tr}\biggl({\mathbb E}e^{\lambda Y_{n-1}}{\mathbb E}e^{\lambda X_n}\biggr)\leq 
{\mathbb E}{\rm tr}\Bigl(e^{\lambda Y_{n-1}}\Bigr)
\Bigl\|{\mathbb E}e^{\lambda X_n}\Bigr\|. 
$$
By induction, we conclude that 
$$
{\mathbb E}{\rm tr}(e^{\lambda Y_n})\leq 
{\mathbb E}{\rm tr}\Bigl(e^{\lambda X_1}\Bigr)
\Bigl\|{\mathbb E}e^{\lambda X_2}\Bigr\|\dots 
\Bigl\|{\mathbb E}e^{\lambda X_n}\Bigr\|. 
$$
Since 
$
{\mathbb E}{\rm tr}\Bigl(e^{\lambda X_1}\Bigr)
=
{\rm tr}\Bigl({\mathbb E}e^{\lambda X_1}\Bigr)
\leq m 
\Bigl\|
{\mathbb E}e^{\lambda X}
\Bigr\|,
$
we get 
\begin{equation}
\label{tri_3}
{\mathbb E}{\rm tr}(e^{\lambda Y_n})\leq 
m
\Bigl\|
{\mathbb E}e^{\lambda X}
\Bigr\|^n.
\end{equation}
It remains to bound the norm $\|{\mathbb E}e^{\lambda X}\|.$ To this end, 
we use Taylor expansion and the condition ${\mathbb E}X=0$ to get 
$$
{\mathbb E}e^{\lambda X}= I_m + 
{\mathbb E}\lambda^2 X^2\biggl[
\frac{1}{2!}+ \frac{\lambda X}{3!}+\frac{\lambda^2 X^2}{4!}+ \dots 
\biggr]\leq
$$
$$ 
I_m+\lambda^2{\mathbb E}X^2\biggl[
\frac{1}{2!}+ \frac{\lambda \|X\|}{3!}+\frac{\lambda^2 \|X\|^2}{4!}+ \dots 
\biggr]=
I_m+\lambda^2{\mathbb E}X^2\biggl[\frac{e^{\lambda \|X\|}-1-\lambda \|X\|}{\lambda ^2 \|X\|^2}
\biggr].
$$
Therefore, for all $\tau>0,$
$$
\Bigl\|{\mathbb E}e^{\lambda X}\Bigr\|\leq 
1+\lambda^2\biggl\|{\mathbb E}X^2\biggl[\frac{e^{\lambda \|X\|}-1-\lambda \|X\|}{\lambda ^2 \|X\|^2}
\biggr]\biggr\|\leq 
$$
$$
1+\lambda^2\Bigl\|{\mathbb E}X^2\Bigr\|
\biggl[\frac{e^{\lambda \tau}-1-\lambda \tau}{\lambda ^2 \tau^2}
\biggr]+
\lambda^2
{\mathbb E}\|X\|^2\biggl[\frac{e^{\lambda \|X\|}-1-\lambda \|X\|}{\lambda ^2 \|X\|^2}
\biggr]I(\|X\|\geq \tau).
$$
Let $M:=2(\log 2)^{1/\alpha}U_X^{(\alpha)}$ and assume that 
$\lambda \leq 1/M.$ Then 
$$
{\mathbb E}\|X\|^2\biggl[
\frac{e^{\lambda \|X\|}-1-\lambda \|X\|}{\lambda ^2 \|X\|^2}
\biggr]I(\|X\|\geq \tau)\leq 
M^2 {\mathbb E}e^{\|X\|/M}I(\|X\|\geq \tau)\leq 
$$
$$
M^2 {\mathbb E}^{1/2}e^{2\|X\|/M}{\mathbb P}^{1/2}\{\|X\|\geq \tau\}.
$$
Since, for $\alpha\geq 1,$ 
$
M=2(\log 2)^{1/\alpha}\Bigl\|\|X\|\Bigr\|_{\psi_{\alpha}}\geq 
2\Bigl\|\|X\|\Bigr\|_{\psi_1}
$
(see van der Vaart and Wellner (1996), p. 95), we have 
${\mathbb E}e^{2\|X\|/M}\leq 2$ and also 
$$
{\mathbb P}\{\|X\|\geq \tau\}\leq 
\exp\biggl\{-2^{\alpha}\log 2\biggl(\frac{\tau}{M}\biggr)^{\alpha}\biggr\}.
$$
As a result, we get the following bound 
$$
\Bigl\|{\mathbb E}e^{\lambda X}\Bigr\|\leq 
1+\lambda^2 \sigma_X^2
\biggl[\frac{e^{\lambda \tau}-1-\lambda \tau}{\lambda ^2 \tau^2}
\biggr]+
2^{1/2}\lambda^2 M^2 \exp\biggl\{-2^{\alpha-1}\log 2\biggl(\frac{\tau}{M}\biggr)^{\alpha}\biggr\}.
$$
Let 
$\tau := M\frac{2^{1/\alpha-1}}{(\log 2)^{1/\alpha}}\log^{1/\alpha} {\frac{M^2}{\sigma_X^2}}$ 
and suppose 
that $\lambda$ satisfies the condition $\lambda \tau \leq 1.$ 
Then, the following bound holds with some constant $C_1>0:$
$$
\Bigl\|{\mathbb E}e^{\lambda X}\Bigr\|\leq 
1+C_1 \lambda^2 \sigma_X^2 \leq \exp\{C_1 \lambda^2 \sigma_X^2\}. 
$$
Thus, we proved that there exist constants $C_1, C_2>0$ such 
that, for all $\lambda $ satisfying the condition 
\begin{equation}
\label{four_4}
\lambda \ U_{X}^{(\alpha)}\biggl(\log {\frac{U_X^{(\alpha)}}{\sigma_X}}\biggr)^{1/\alpha}\leq C_2,
\end{equation}
we have
$
\Bigl\|{\mathbb E}e^{\lambda X}\Bigr\|\leq \exp\{C_1 \lambda^2 \sigma_X^2\}. 
$
This can be combined with (\ref{odin_1}), (\ref{dwa_2}) and (\ref{tri_3})
to get 
$$
{\mathbb P}\{\|Y_n\|\geq t\}\leq 2m \exp\Bigl\{-\lambda t + C_1 \lambda^2 n \sigma_X^2\Bigr\}.
$$
It remains now to minimize the last bound with respect to all $\lambda$
satisfying (\ref{four_4}) to get that, for some constant $K>0,$ 
$$
{\mathbb P}\{\|Y_n\|\geq t\}\leq 2m \exp\biggl\{-\frac{1}{K}
\frac{t^2}{n\sigma_X^2+ t U_X^{(\alpha)} \log^{1/\alpha}(U_X^{(\alpha)}/\sigma_X)}\biggr\},
$$
which immediately implies (\ref{bernstein-2}).

\qed

Note that, in a standard way, one can deduce bounds 
on the expectation from the exponential bounds on tail probabilities.
In particular, (\ref{bernstein-0}) implies that 
 \begin{equation}
\label{bernstein-1_A}
{\mathbb E}\biggl\|\frac{X_1+\dots+X_n}{n}\biggr\|\leq 
C\biggl(\sigma_X\sqrt{\frac{\log(2m)}{n}}\bigvee U\frac{\log(2m)}{n}\biggr). 
\end{equation}
Similarly, Proposition \ref{bernstein_A} implies that 
\begin{equation}
\label{bernstein-2_A}
{\mathbb E}\biggl\|\frac{X_1+\dots+X_n}{n}\biggr\|\leq 
C\biggl(\sigma_X\sqrt{\frac{\log(2m)}{n}}\bigvee 
U_{X}^{(\alpha)}\biggl(\log \frac{U_X^{(\alpha)}}{\sigma_X}\biggr)^{1/\alpha}\frac{\log(2m)}{n}\biggr) 
\end{equation}
Combining the last bounds with Talagrand's concentration inequality 
leads to somewhat different versions of bounds (\ref{bernstein-1})
and (\ref{bernstein-2}) that can be better in some applications. 
Namely, denote 
$$
\tilde \sigma_X^2:=\sup_{u,v\in {\mathbb C}^m,|u|\leq 1, |v|\leq 1}
{\mathbb E}|\langle Xu,v\rangle|^2.
$$
It is easy to check that 
$
\tilde \sigma_X^2 \leq \sigma_X^2.
$
Moreover, in some cases, it can be significantly smaller (for instance,
if $X$ is sampled at random from the matrix completion basis, then $\sigma_X^2$ is of the order $m^{-1}$ and $\tilde \sigma_X^2$ 
is equal to $m^{-2}$). The expectation bound (\ref{bernstein-1_A})
and Talagrand's concentration inequality imply that with probability 
at least $1-e^{-t}$ 
\begin{equation}
\label{bernstein-1_B}
\biggl\|\frac{X_1+\dots+X_n}{n}\biggr\|
\leq 
C\biggl(
\sigma_X\sqrt{\frac{\log(2m)}{n}}\bigvee 
\tilde \sigma_X\sqrt{\frac{t}{n}}\bigvee U\frac{\log(2m)}{n}
\bigvee U\frac{t}{n}
\biggr). 
\end{equation}
Similarly, combining the expectation bound (\ref{bernstein-2_A})
for $\alpha=1$ with Adamczak's version of Talagrand's inequality 
(see Section 3.2), we get that with probability at least $1-e^{-t}$
\begin{equation}
\label{bernstein-2_B}
\biggl\|\frac{X_1+\dots+X_n}{n}\biggr\|\leq 
C\biggl(\sigma_X\sqrt{\frac{\log(2m)}{n}}\bigvee \tilde \sigma_X\sqrt{\frac{t}{n}}\bigvee 
U_{X}^{(1)}\biggl(\log \frac{U_X^{(1)}}{\sigma_X}\biggr)\frac{\log(2m)}{n}
\bigvee U_X^{(1)}\frac{t\log n}{n}
\biggr).
\end{equation}
In the examples when $\tilde \sigma_X^2$ is significantly smaller than 
$\sigma_X^2,$ these bounds might be better than (\ref{bernstein-1}) and 
(\ref{bernstein-2}), especially, when they are used for large values 
of $t.$

In principle, using bounds (\ref{bernstein-1_B}) and (\ref{bernstein-2_B})
in the proofs of the following sections instead of (\ref{bernstein-1}) and (\ref{bernstein-2}) provides a way  
to obtain probabilistic oracle inequalities with probabilities of the 
error decreasing exponentially with $m $ or $n$ (this is the way in which 
error bounds are written in the papers by Candes and Plan (2009) and Rohde and Tsybakov (2009)). We are not pursuing this approach here.

\section{Approximation Error}

A natural first step in the analysis of the problem is to study its version with the true risk instead of the empirical risk. 
The true risk with respect to the quadratic loss is equal to 
$$
{\mathbb E}(Y-\langle S,X\rangle)^2= {\mathbb E}(\langle \rho, X\rangle +\xi -\langle S,X\rangle)^2={\mathbb E}\langle S-\rho, X \rangle^2 + {\mathbb E}\xi^2,
$$
where we used the assumptions that $X$ and $\xi$ are independent 
and ${\mathbb E}\xi=0.$ Thus, the penalized true risk minimization problem 
becomes 
\begin{equation}
\label{entropy_penalty_true}
\rho^{\eps}:={\rm argmin}_{S\in {\cal S}}\biggl[
{\mathbb E}\langle S-\rho, X \rangle^2 + \eps\ {\rm tr}(S\log S)
\biggr]
\end{equation}
and the goal is to study the error of approximation of $\rho$ by $\rho^{\eps}$ depending on the value of regularization parameter 
$\eps>0.$ 
The next propositions show that if there exists an oracle 
$S\in {\cal S}$ that provides a good approximation of the true state 
$\rho $ in a sense that $\|S-\rho\|_{L_2(\Pi)}$ is small, then $\rho^{\eps}$ belongs to an $L_2(\Pi)$-ball around $S$ of small enough radius that can be controlled in terms of 
the operator norm $\|\log S\|$ or in terms of more subtle characteristics
of the oracle $S.$
They also provide upper bounds  
on the approximation error $\|\rho^{\eps}-\rho\|_{L_2(\Pi)}.$

We will first obtain a simple bound on  
$\|\rho^{\eps}-S\|_{L_2(\Pi)}$ 
for an arbitrary oracle $S\in {\cal S}$ of full rank 
expressed in terms of the operator norm $\|\log S\|$ of its logarithm.
For simplicity, we assume that $\|\log S\|=+\infty$ in the case 
when ${\rm rank}(S)<m$ (and $\log S$ is not defined). Note, however,
that ${\rm tr}(S\log S)$ is well defined and finite even in the case
when ${\rm rank}(S)<m.$

\begin{proposition}
\label{prop3}
For all $S\in {\cal S},$
$
\|\rho^{\eps}-S\|_{L_2(\Pi)}\leq \|S-\rho\|_{L_2(\Pi)}
+\sqrt{\eps \|\log S\|}.
$
This implies that 
$$
\|\rho^{\eps}-\rho\|_{L_2(\Pi)}\leq 2\|S-\rho\|_{L_2(\Pi)}
+\sqrt{\eps \|\log S\|},
$$
and, in particular, for $S=\rho,$
$
\|\rho^{\eps}-\rho\|_{L_2(\Pi)}^2\leq 
\eps \|\log \rho\|.
$
\end{proposition}

For a differentiable mapping $g$ from an open subset $G\subset {\mathbb M}_m({\mathbb C})$ into ${\mathbb M}_m({\mathbb C}),$ denote by 
$Dg(A;H)$ its differential at a matrix $A\in G$ in the direction $H\in {\mathbb M}_m({\mathbb C}),$
that is, $Dg(A;H)$ is linear with respect to $H$ and  
$$
g(A+H)=g(A)+ Dg(A;H)+ o(\|H\|)\ {\rm as}\ \|H\|\to 0.
$$ 
The following lemma is a simple 
corollary of Theorem V.3.3 in Bhatia (1996):

\begin{lemma}
\label{from_Bhatia}
Let $f$ be a function continuously differentiable in an open interval 
$I\subset {\mathbb R}.$ Suppose that $A$ is a Hermitian matrix 
whose spectrum belongs to $I.$ Then the mapping $B\mapsto g(B):={\rm tr}(f(B))$ 
is differentiable at $A$ and 
$
Dg(A;H)={\rm tr}(f^{\prime}(A)H).
$
\end{lemma}

{\bf Proof of Proposition \ref{prop3}.} 
Denote the penalized risk 
$$
L(S):={\mathbb E}\langle S-\rho, X \rangle^2 + \eps\ {\rm tr}(S\log S).
$$
It is easy to see that the solution $\rho^{\eps}$ of 
problem (\ref{entropy_penalty_true}) is a full rank matrix.
To prove this, assume that ${\rm rank}(\rho^{\eps})<m.$
Let $\tilde \rho :=(1-\delta)\rho^{\eps}+\delta I_m,$ where $I_m$ is the $m\times m$ identity matrix. Then, for small enough $\delta,$ $\tilde \rho$ is a full rank matrix and it is straightforward to show that 
the penalized risk $L(\tilde \rho)$ is strictly smaller than $L(\rho^{\eps})$ (for some small $\delta>0$). 
It is also easy to check that, for
any $S\in {\cal S}$ of full rank, $\log S$ is well defined and 
the differential of the functional $L$ in the 
direction $\nu\in {\mathbb M}_m({\mathbb C})$ is equal to  
$$
DL(S;\nu) = 
2{\mathbb E}\langle S-\rho,X\rangle \langle \nu, X\rangle + \eps\ {\rm tr}(\nu \log S).
$$
This follows from the fact that the first term of the functional $L$
is differentiable since it is quadratic. The differentiability of the 
penalty term is based on Lemma \ref{from_Bhatia} 
(it is enough to apply 
this lemma to the function $f(u)=u\log u$). 
Since $\rho^{\eps}$ is the minimal point of $L$ in ${\cal S},$ 
we can conclude that, for an arbitrary $S\in {\cal S},$
$
DL(\rho^{\eps};S-\rho^{\eps})\geq 0.
$
This implies that 
$$
DL(S;S-\rho^{\eps})-DL(\rho^{\eps};S-\rho^{\eps})\leq DL(S;S-\rho^{\eps}),
$$
which, by a simple algebra, becomes 
\begin{equation}
\label{BB1}
2\|S-\rho^{\eps}\|_{L_2(\Pi)}^2+ \eps K(S;\rho^{\eps})
\leq 
2{\mathbb E}\langle S-\rho,X\rangle \langle S-\rho^{\eps}, X\rangle + \eps\ \langle S-\rho^{\eps},\log S \rangle.
\end{equation}
To conclude the proof, note that (\ref{BB1}), the bound 
$\|S-\rho^{\eps}\|_1\leq 2$ and Cauchy-Schwarz inequality imply that
$$
2\|S-\rho^{\eps}\|_{L_2(\Pi)}^2+ \eps K(S;\rho^{\eps})
\leq 
2\|S-\rho\|_{L_2(\Pi)}\|S-\rho^{\eps}\|_{L_2(\Pi)} + 2\eps \|\log S\|.
$$
Solving the last inequality with respect to $\|\rho^{\eps}-S\|_{L_2(\Pi)}$
and using the fact that $K(S;\rho^{\eps})\geq 0,$
yields the bound 
$$
\|\rho^{\eps}-S\|_{L_2(\Pi)}\leq \frac{\|S-\rho\|_{L_2(\Pi)}}{2}
+\sqrt{\frac{\|S-\rho\|_{L_2(\Pi)}^2}{4}+\eps \|\log S\|},
$$
which implies 
$
\|\rho^{\eps}-S\|_{L_2(\Pi)}\leq \|S-\rho\|_{L_2(\Pi)}
+\sqrt{\eps \|\log S\|},
$
and the result follows.

\qed

To obtain more subtle bounds with approximation error of the order 
$O(\eps^2)$ instead of $O(\eps),$  
we introduce and use the following quantity
$$
a(W):=a_{\Pi}(W):=a_X(W):=\sup\biggl\{\langle W, U\rangle : U\in {\mathbb M}_{m}({\mathbb C}),
U=U^{\ast}, {\rm tr}(U)=0, 
\|U\|_{L_2(\Pi)}=1
\biggr\},
$$
which will be called the {\it alignment coefficient} of $W.$
Similar quantities were used in the commutative case (Koltchinskii (2009)).
Note that, for all constants $c,$
\begin{equation}
\label{shift}
a(W+cI_m)=a(W)
\end{equation}
(since $\langle I_m, U\rangle=0$ for all $U$ of zero trace).
In addition, we have 
\begin{equation}
\label{scale}
a_{cX}(W)=\frac{1}{|c|}a_X(W),\ c\neq 0. 
\end{equation}

Let $\{E_i:i=1,\dots, m^2\}$ be an orthonormal 
basis of ${\mathbb M}_m({\mathbb C})$ consisting of Hermitian 
matrices and let ${\cal K}:=\Bigl(\langle E_j, E_k\rangle_{L_2(\Pi)}\Bigr)_{j,k=1}^{m^2}$ be the Gram 
matrix of the functions $\{\langle E_j,\cdot\rangle : j=1,\dots, m^2\}$
in the space $L_2(\Pi).$ Clearly, the mapping 
$J:{\mathbb M}_m({\mathbb C})\mapsto \ell_2^{m^2}({\mathbb C}),$
$$
J U = \Bigl(\langle U,E_j\rangle: j=1,\dots, m^2\Bigr),\ U\in {\mathbb M}_m({\mathbb C}),
$$
is an isometry. If now we define 
$\bar {\cal K}:{\mathbb M}_m({\mathbb C})\mapsto {\mathbb M}_m({\mathbb C})$ as $\bar {\cal K}:=J^{-1}{\cal K}J,$ then we also 
have 
$\bar {\cal K}^{1/2}=J^{-1}{\cal K}^{1/2}J,\ \bar {\cal K}^{-1/2}=J^{-1}{\cal K}^{-1/2}J.$ 
As a consequence, for any matrix $U=\sum_{j=1}^{m^2}u_j E_j,$
$$
\|U\|_{L_2(\Pi)}^2 = \sum_{j,k=1}^{m^2}\langle E_j, E_k\rangle_{L_2(\Pi)}
u_j \bar u_k=\langle {\cal K}u,u\rangle_{\ell_2}=
\|{\cal K}^{1/2} u\|_{\ell_2}^2=\|\bar {\cal K}^{1/2}U\|_2^2,
$$
and it is not hard to conclude that 
$
a(W)\leq \|\bar {\cal K}^{-1/2} W\|_{2}.
$
Moreover, in view of (\ref{shift}), for an arbitrary 
scalar $c,$
$$
a(W)\leq \|\bar {\cal K}^{-1/2} (W+cI_m)\|_{2}.
$$
This shows that the size of $a(W)$ depends on how $W$ is ``aligned''
with the eigenspaces of the Gram matrix ${\cal K}.$ In a special case 
when, for all $A,$ $\|A\|_{L_2(\Pi)}=\|A\|_2,$
the functions $\{\langle E_j,\cdot\rangle : j=1,\dots, m^2\}$
form an orthonormal system in the space $L_2(\Pi)$ and 
the Gram matrix ${\cal K}$ is the identity matrix. 
In this case, we simply have the bound 
$$
a(W)\leq \inf_{c}\|W+cI_m\|_{2}.
$$

In the next statement, we use the alignment coefficient $a(\log S)$ to control 
the $L_2(\Pi)$-distance $\|\rho^{\eps}-S\|_{L_2(\Pi)}$ and 
the Kullback-Leibler ``distance'' $K(\rho^{\eps};S)$ from 
the true solution $\rho^{\eps}$ to an arbitrary oracle $S.$

\begin{proposition}
\label{prop1}
For all $S\in {\cal S},$
$$
\|\rho^{\eps}-S\|_{L_2(\Pi)}^2+ 
\frac{\eps}{2} K(\rho^{\eps};S)
\leq 
\biggl(\|S-\rho\|_{L_2(\Pi)} + \frac{\eps}{2} a(\log S)\biggr)^2.
$$
In particular, it implies that 
$
\|\rho^{\eps}-\rho\|_{L_2(\Pi)}^2 + 
\frac{\eps}{2} K(\rho^{\eps};\rho)
\leq 
\frac{\eps^2}{4} a^2(\log \rho).
$
Moreover, the following bound also holds:
$$
\|\rho^{\eps}-\rho\|_{L_2(\Pi)}^2
\leq 
\inf_{S\in {\cal S}}\biggl[\|S-\rho\|_{L_2(\Pi)}^2 + 
\eps a(\log S)\|S-\rho\|_{L_2(\Pi)}+
\frac{\eps^2}{2} a^2(\log S)\biggr].
$$
\end{proposition}

{\bf Proof}. Our starting point is the relationship (\ref{BB1})
from the proof of Proposition \ref{prop3}. 
It follows from the definition of $a(W),$ from (\ref{BB1}) 
and from Cauchy-Schwarz inequality that 
$$
2\|S-\rho^{\eps}\|_{L_2(\Pi)}^2+ \eps K(S;\rho^{\eps})
\leq 
2\|S-\rho\|_{L_2(\Pi)}\|S-\rho^{\eps}\|_{L_2(\Pi)} + \eps a(\log S)\|S-\rho^{\eps}\|_{L_2(\Pi)}.
$$
It remains to solve the last inequality for $\|S-\rho^{\eps}\|_{L_2(\Pi)}$
to obtain the first bound of the proposition. 
The second bound is its special case with $S=\rho.$
To prove the third bound note that, by the definition of $\rho^{\eps},$ for all $S\in {\cal S},$
$$
\|\rho^{\eps}-\rho\|_{L_2(\Pi)}^2 + \eps {\rm tr}(\rho^{\eps}\log \rho^{\eps})
\leq \|S-\rho\|_{L_2(\Pi)}^2 + \eps {\rm tr}(S\log S),
$$
which implies 
$$
\|\rho^{\eps}-\rho\|_{L_2(\Pi)}^2  
\leq 
\|S-\rho\|_{L_2(\Pi)}^2 + \eps 
({\rm tr}(S\log S)-{\rm tr}(\rho^{\eps}\log \rho^{\eps}))\leq 
$$
$$
\|S-\rho\|_{L_2(\Pi)}^2 + 
\eps {\rm tr}(\log S (S-\rho^{\eps}))\leq 
\|S-\rho\|_{L_2(\Pi)}^2 + 
\eps a(\log S)
\|\rho^{\eps}-S\|_{L_2(\Pi)},
$$
where we used the fact that, by convexity of the function
$S\mapsto {\rm tr}(S\log S),$ 
$$
{\rm tr}(S\log S)-{\rm tr}(\rho^{\eps}\log \rho^{\eps})
\leq {\rm tr}(\log S (S-\rho^{\eps})).
$$
It remains to bound 
$\|\rho^{\eps}-S\|_{L_2(\Pi)}$ 
from above using the first inequality of the proposition.  

\qed

A consequence of propositions \ref{prop3} and \ref{prop1} is 
that 
\begin{equation}
\label{ooo}
\|\rho^{\eps}-\rho\|_{L_2(\Pi)}^2 \leq 
\frac{\eps^2}{4} a^2(\log \rho)\wedge \eps \|\log \rho\|.
\end{equation}

We will now provide versions of approximation error 
bounds for special types of oracles $S\in {\cal S}.$

{\bf Low Rank Oracles}. First we show how to adapt the bounds of Proposition \ref{prop1} expressed in terms of the alignment coefficient $a(\log S)$ for a full rank matrix $S$ (for which $\log S$ is well defined)
to the case when $S$ is an oracle of a small rank $r<m.$
For a subspace $L$ of ${\mathbb C}^m,$
denote 
$
\Lambda (L):= \sup_{\|A\|_{L_2(\Pi)}\leq 1}\|P_L A P_L\|_{2}.
$
Suppose that $S\in {\cal S}$ 
is a matrix of rank $r.$ To be specific, let 
$
S = \sum_{j=1}^r \gamma_j (e_j\otimes e_j),
$
where $\gamma_j$ are positive eigenvalues of $S$ and 
$\{e_1,\dots, e_m\}$ is an orthonormal basis of ${\mathbb C}^m.$
Let $L$ be the linear span of the vectors $e_1,\dots, e_r.$

\begin{proposition}
\label{low_rank_approx}
There exists a numerical constant $C>0$ such that, for all $\eps>0,$
$$
\|\rho^{\eps}-\rho\|_{L_2(\Pi)}^2
\leq 2\|S-\rho\|_{L_2(\Pi)}^2+
C\eps^2\biggl[\Lambda^2(L) r  
\log^2 \biggl(1+\frac{m}{\eps\wedge 1}\biggr)+
{\mathbb E}\|X\|^2\biggr].
$$
\end{proposition}

{\bf Proof}. 
Note that, for all matrices $W$ of rank $r$ 
``supported'' in the space $L$
in the sense that $W=P_L W P_L,$ we have  
$$
a(W)\leq 
\sup_{\|U\|_{L_2(\Pi)}\leq 1} \langle W,U\rangle
=
\sup_{\|U\|_{L_2(\Pi)}\leq 1} \langle W, P_L U P_L\rangle
\leq  \Lambda (L)\|W\|_2. 
$$
For $\delta \in (0,1),$ consider 
$
S_{\delta}:=(1-\delta)S + \delta \frac{I_m}{m}.
$ 
Then, using the fact that $a(W+cI_m)=a(W),$ we get 
$$
\log S_{\delta} = 
\sum_{j=1}^r \Bigl(\log((1-\delta)\gamma_j+\delta/m)-\log (\delta/m)\Bigr)(e_j\otimes e_j)
+ \log (\delta/m)I_m 
$$ 
and 
$$
a(\log S_{\delta})= 
a\biggl(\sum_{j=1}^r \Bigl(\log((1-\delta)\gamma_j+\delta/m)-\log (\delta/m)\Bigr)(e_j\otimes e_j)\biggr)\leq 
$$
$$
\Lambda(L) \biggl\|\sum_{j=1}^r \Bigl(\log((1-\delta)\gamma_j+\delta/m)-\log (\delta/m)\Bigr)(e_j\otimes e_j)\biggr\|_2\leq 
$$
$$
\Lambda (L)\biggl(\sum_{j=1}^r \log^2 \biggl(1+ \frac{m\gamma_j}{\delta}\biggr)\biggr)^{1/2}
\leq \Lambda (L)\sqrt{r}\log \biggl(1+\frac{m\|S\|}{\delta}\biggr).
$$
Note also that 
$
\|S-S_{\delta}\|_{L_2(\Pi)}^2 = \delta^2 \|S -I_m/m\|_{L_2(\Pi)}^2 
\leq 4 \delta^2 {\mathbb E}\|X\|^2,
$
since 
$$\|S -I_m/m\|_{L_2(\Pi)}^2\leq 
2({\mathbb E}\langle S,X\rangle^2 +{\mathbb E}\langle I_m/m,X\rangle^2)\leq 
$$
$$
2(\|S\|_1^2{\mathbb E}\|X\|^2 + \|I_m/m\|_1^2{\mathbb E}\|X\|^2)\leq 
4 {\mathbb E}\|X\|^2.
$$ 
Thus, it easily follows from Proposition \ref{prop1} that 
$$
\|\rho^{\eps}-\rho\|_{L_2(\Pi)}^2
\leq  
\frac{3}{2}\|S_{\delta}-\rho\|_{L_2(\Pi)}^{2}+ \eps^2 a^2(\log S_{\delta})
\leq
$$
$$
\frac{3}{2}\Bigl(\|S-\rho\|_{L_2(\Pi)}+\|S_{\delta}-S\|_{L_2(\Pi)}\Bigr)^2+
\Lambda^2(L) r\eps^2
\log^2 \biggl(1+\frac{m}{\delta}\biggr)
\leq
$$
$$
\frac{3}{2}
\Bigl(\frac{4}{3}\|S-\rho\|_{L_2(\Pi)}^2+4\|S_{\delta}-S\|_{L_2(\Pi)}^2
\Bigr)+
\Lambda^2(L) r\eps^2
\log^2 \biggl(1+\frac{m}{\delta}\biggr)
\leq
$$
$$
2\|S-\rho\|_{L_2(\Pi)}^2+24 {\mathbb E}\|X\|^2\delta^2 + 
\Lambda^2(L) r\eps^2
\log^2 \biggl(1+\frac{m}{\delta}\biggr).
$$
Taking $\delta =\eps \wedge 1,$ this yields 
$$
\|\rho^{\eps}-\rho\|_{L_2(\Pi)}^2
\leq 2\|S-\rho\|_{L_2(\Pi)}^2+
C\eps^2\biggl[\Lambda^2(L)r  
\log^2 \biggl(1+
\frac{m}{\eps\wedge 1}
\biggr)+{\mathbb E}\|X\|^2\biggr]
$$
with a numerical constant $C>0.$

\qed

{\bf Remark}. The bound of Proposition \ref{low_rank_approx} can be 
also written in the following form that might be preferable when 
${\mathbb E}\|X\|^2$ is large:
$$
\|\rho^{\eps}-\rho\|_{L_2(\Pi)}^2
\leq 2\|S-\rho\|_{L_2(\Pi)}^2+
C\eps^2\biggl[\Lambda^2(L) r  
\log^2 \biggl(1+\biggl(\frac{m{\mathbb E}^{1/2}\|X\|^2}{\eps}\vee m\biggr)\biggr)+1\biggr].
$$
In the proof, it is enough to take 
$
\delta:=\frac{\eps}{{\mathbb E}^{1/2}\|X\|^2}\wedge 1.
$

Note that if 
$\{E_i,\ i=1,\dots, m^2\}$ 
is an orthonormal basis 
of ${\mathbb M}_m({\mathbb C})$ consisting of Hermitian matrices
and $X$ is uniformly distributed in  
$\{E_i,\ i=1,\dots, m^2\},$ then for all Hermitian $A$
$$
\|A\|_{L_2(\Pi)}^2 = {\mathbb E}\langle A,X\rangle^2 
=m^{-2}\sum_{j=1}^{m^2} \langle A,E_j\rangle^2=
m^{-2} \|A\|_2^2.
$$
Therefore 
$
\Lambda (L)\leq 
\sup_{\|A\|_{L_2(\Pi)}\leq 1}\|A\|_2 =
\sup_{\|A\|_{2}\leq m}\|A\|_2 =m.
$
Also, in this case $\|X\|\leq \|X\|_2=1.$
Thus, Proposition \ref{low_rank_approx} yields 
$$
\|\rho^{\eps}-\rho\|_{L_2(\Pi)}^2
\leq 2\|S-\rho\|_{L_2(\Pi)}^2+
C m^2 r \eps^2 
\log^2 \biggl(1+\frac{m}{\eps \wedge 1}\biggr)+C\eps^2.
$$

{\bf Gibbs Oracles}. Let $H$ be a Hermitian 
matrix (``a Hamiltonian'') and let $\beta>0.$ Consider the following 
density matrix (a ``Gibbs oracle''): 
$$
\rho_{H,\beta}:= \frac{e^{-\beta H}}{{\rm tr}(e^{-\beta H})}.
$$
For simplicity, assume in what follows that $\beta=1$ (in fact, one can always replace $H$ by $\beta H$) and denote 
$
\rho_H:=\frac{e^{-H}}{{\rm tr}(e^{-H})}.
$ 
Let $\gamma_1\leq \gamma_2 \leq \dots \leq \gamma_m$
be the eigenvalues of $H$ and $e_1,\dots, e_m$ be the corresponding 
eigenvectors. Let $L_r={\rm l.s.}(\{e_1,\dots, e_r\})$ and 
$
H_{\leq r}:= \sum_{j=1}^r \gamma_j (e_j\otimes e_j),\ \ 
H_{> r}:= \sum_{j=r+1}^m \gamma_j (e_j\otimes e_j).
$ 
It is easy to see that 
$$
\|P_{L_{r}^{\perp}} \rho_{H} P_{L_{r}^{\perp}}\|_1 = \frac{\sum_{k\geq r+1}e^{-\gamma_k}}{\sum_{k\geq 1}e^{-\gamma_k}}=:
\delta_r (H).
$$ 
Under reasonable conditions on the spectrum of $H,$ the quantity $\delta_r(H)$ decreases 
fast enough when $r$ increases. Thus, $\rho_{H}$ can be well approximated by low rank matrices.

The next statement follows immediately from Proposition 
\ref{prop1}. Here the unknown density matrix $\rho$ is approximated 
by a Gibbs model with an arbitrary Hamiltonian. The error is controlled in terms of the $L_2(\Pi)$-distance between $\rho$ and the oracle $\rho_{H}$
and also in terms of the alignment coefficient $a(H_{\leq r})$ for 
a ``low rank part'' $H_{\leq r}$ of the Hamiltonian $H$ and 
the quantity $\delta_r(H).$

\begin{proposition}
\label{Gibbs_approx}
For all Hermitian 
nonnegatively definite matrices $H$ and 
for all $\eps>0,$
$$
\|\rho^{\eps}-\rho\|_{L_2(\Pi)}^2
\leq 2\|\rho_{H}-\rho\|_{L_2(\Pi)}^2+ 24
\max_{1\leq k\leq m}{\mathbb E}\langle X e_k, e_k\rangle^2
\delta_r^2 (H)+
a^2(H_{\leq r}) \eps^2.
$$
\end{proposition}

{\bf Proof}. We will use the last bound of proposition \ref{prop1}
with $S=\rho_{H_{\leq r}}.$ Note that 
$$
a(\log \rho_{H_{\leq r}})= a(-H_{\leq r} - \log {\rm tr}(e^{-H_{\leq r}})I_m)= a(H_{\leq r}).
$$
Therefore, we have 
$$
\|\rho^{\eps}-\rho\|_{L_2(\Pi)}^2
\leq 
\|\rho_{H_{\leq r}}-\rho\|_{L_2(\Pi)}^2 + 
\eps a(H_{\leq r})\|\rho_{H_{\leq r}}-\rho\|_{L_2(\Pi)}+
\frac{\eps^2}{2} a^2(H_{\leq r})\leq 
$$
$$
\frac{3}{2}\|\rho_{H_{\leq r}}-\rho\|_{L_2(\Pi)}^2 + 
\eps^2 a^2(H_{\leq r}).
$$
In addition to this,
$$
\|\rho_{H}-\rho_{H_{\leq r}}\|_{L_2(\Pi)}
=\biggl\|
\frac{\sum_{k=1}^m e^{-\gamma_k}(e_k\otimes e_k)}{\sum_{k=1}^m e^{- \gamma_k}} -
\frac{\sum_{k=1}^r e^{-\gamma_k}(e_k\otimes e_k)}{\sum_{k=1}^r e^{- \gamma_k}} 
\biggr\|_{L_2(\Pi)},
$$
which can be easily bounded from above by 
$$
2 \delta_r(H)\max_{1\leq k\leq m}\|e_k\otimes e_k\|_{L_2(\Pi)}=
2 \delta_r(H)\max_{1\leq k\leq m}
{\mathbb E}^{1/2}\langle X e_k, e_k\rangle^2.
$$
The result follows immediately (by the same argument as in the proof 
of Proposition \ref{low_rank_approx}).

\qed



\section{Random Error Bounds and Oracle Inequalities}

We now turn to the analysis of random error of the estimator
$\hat \rho^{\eps}.$ We obtain upper bounds on the $L_2(\Pi)$
and Kullback-Leibler distances of this estimator to an arbitrary 
oracle $S\in {\cal S}$ of full rank. 
In particular, this includes bounding the distances between $\hat \rho^{\eps}$ and $\rho^{\eps}.$ 
As a consequence, we will obtain {\bf oracle inequalities} for the 
empirical solution $\hat \rho^{\eps}.$  
The size 
of both errors $\|\hat \rho^{\eps}-S\|_{L_2(\Pi)}^2$ and 
$K(\hat \rho^{\eps};S)$ will be controlled in terms of 
the squared $L_2(\Pi)$-distance $\|S-\rho\|_{L_2(\Pi)}^2$
from the oracle to the target density matrix $\rho$ and also in terms 
of such characteristics of the oracle as the norm $\|\log S\|$
or the alignment coefficient $a(\log S)$ that have been already
used in the approximation error bounds of the previous section
(see propositions \ref{prop3}, \ref{prop1}). However, in the case
of the random error, we also need some additional quantities that 
describe the properties of the design distribution $\Pi$ and 
of the noise $\xi.$  These quantities are explicitly involved in the statements of the results below which makes these statements somewhat  
complicated. At the same time, it is easy to control these quantities in concrete examples and to derive in special 
cases the bounds that are easier to understand.

{\bf Assumptions on the design distribution $\Pi$}. 
In this section, it will be assumed that
$X$ is a random Hermitian $m\times m$  matrix 
and that, for some constant $U>0,$ $\|X\|\leq U.$ We will denote 
$$
\sigma_X^2:=\|{\mathbb E}(X-{\mathbb E}X)^2\|, \ \  
\sigma_{X\otimes X}^2:=\|{\mathbb E}(X\otimes X-{\mathbb E}(X\otimes X))^2\|.
$$
Let $L\subset {\mathbb C}^m$ be a subspace of dimension $r\leq m $
and let 
${\cal P}_L:{\mathbb M}_m({\mathbb C})\mapsto {\mathbb M}_m({\mathbb C}),$
$
{\cal P}_L x := x-P_{L^{\perp}} x P_{L^{\perp}}.
$
We will use the following quantity:
$$
\beta (L):= \sup_{\|A\|_{L_2(\Pi)}\leq 1}\|{\cal P}_L A\|_{L_2(\Pi)}.
$$
Note that $\|{\cal P}_L A\|_2\leq \|A\|_2$ (for a proof, choose a basis 
$\{e_1,\dots, e_m\}$ of ${\mathbb C}^m$ such that $L={\rm l.s.}(e_1,\dots, e_r)$ and represent the linear transformations in this basis). 
If, for all $A,$ 
$
K_1\|A\|_2 \leq \|A\|_{L_2(\Pi)}\leq K_2 \|A\|_2,
$
then $\beta (L)\leq K_2/K_1.$ In particular, if $K_1=K_2,$ then 
$\beta (L)=1$ (which is the case, for instance, when $X$ is sampled at random from an orthonormal basis).

{\bf Assumptions on the noise $\xi$}. 
Recall that ${\mathbb E}\xi=0$ and let 
$
\sigma_{\xi}^2 := {\mathbb E}\xi^2<+\infty.
$
We will further assume that the noise is uniformly 
bounded by a constant $c_{\xi}>0:$ $|\xi|\leq c_{\xi},$
and the proofs of the results of this section will be given 
under this assumption. Alternatively, one can assume that 
the noise is not necessarily uniformly bounded,
but $\|\xi\|_{\psi_1}<+\infty.$ This includes, for instance,
the case of Gaussian noise. 
For such an unbounded noise, one should replace
in the proofs of theorems \ref{simple}, \ref{oracle} and \ref{random_error_1} below the noncommutative Bernstein inequality of Ahlswede and Winter by the bound 
of Proposition \ref{bernstein_A}. One should also use 
a version of concentration inequality for empirical processes 
by Adamczak (2008) instead of the usual version of Talagrand
for bounded function classes (see Section 3).

Given $t>0,$ denote 
$t_m:=t+\log (2m), \ \ \tau_n :=t+\log\log_2(2n)
$
and 
$$
\eps_{n,m}:=(\sigma_{\xi} \sigma_X \vee \sigma_{\xi}\|{\mathbb E}X\|\vee \sigma_{X\otimes X})
\sqrt{\frac{t_{m}}{n}}\bigvee
(c_{\xi}U\vee U^2)\frac{t_{m}}{n}. 
$$
We will start with a simple result in spirit of approximation 
error bound of Proposition \ref{prop3}.

\begin{theorem}
\label{simple}
There exists a constant $C>0$ such that, for all $S\in {\cal S}$
and for all $\eps\geq 0,$
with probability at least $1-e^{-t}$
\begin{eqnarray}
\label{simple_1}
&&
\nonumber
\|\hat \rho^{\eps}-S\|_{L_2(\Pi)}^2\leq \|S-\rho\|_{L_2(\Pi)}^2 
+C\biggl[\eps 
(\|\log S\|\bigwedge 
\log \Gamma )
\bigvee 
\|S-\rho\|_{L_2(\Pi)}U\sqrt{\frac{t_m}{n}}\bigvee
\\
&&
(\sigma_{\xi}\sigma_X \vee \sigma_{\xi}
\|{\mathbb E}X\|\vee \sigma_{X\otimes X})\sqrt{\frac{t_m}{n}}\bigvee 
(c_{\xi}U\vee U^2)\frac{t_m}{n} 
\biggr]
\end{eqnarray}
and 
\begin{eqnarray}
\label{simple_2}
&&
\nonumber
\|\hat \rho^{\eps}-\rho\|_{L_2(\Pi)}^2\leq 
\|S-\rho\|_{L_2(\Pi)}^2 
+C\biggl[\eps 
(\|\log S\|\wedge 
\log \Gamma)
\bigvee 
\|S-\rho\|_{L_2(\Pi)}U\sqrt{\frac{t_m}{n}}\bigvee
\\
&&
(\sigma_{\xi}\sigma_X \vee \sigma_{\xi}
\|{\mathbb E}X\|\vee \sigma_{X\otimes X})\sqrt{\frac{t_m}{n}}\bigvee 
(c_{\xi}U\vee U^2)\frac{t_m}{n} 
\biggr],
\end{eqnarray}
where 
$
\Gamma:=\frac{m {\mathbb E}^{1/2}\|X\|^2}{\sqrt{\eps}}\vee m.
$
In particular,
\begin{equation}
\label{simple_3}
\|\hat \rho^{\eps}-\rho\|_{L_2(\Pi)}^2\leq 
C\biggl[\eps 
(\|\log \rho\|\wedge 
\log \Gamma)\bigvee 
(\sigma_{\xi}\sigma_X \vee \sigma_{\xi}
\|{\mathbb E}X\|\vee \sigma_{X\otimes X})\sqrt{\frac{t_m}{n}}\bigvee 
(c_{\xi}U\vee U^2)\frac{t_m}{n} 
\biggr].
\end{equation}
\end{theorem}

Note that this result holds for all $\eps\geq 0,$ including 
the case of $\eps=0$ that corresponds to the least squares 
estimator over the set ${\cal S}$ of all density matrices. 
The approximation error term $\|\log S\|\eps$ in the 
bounds of Theorem \ref{simple} is of the order $O(\eps)$
(as in Proposition \ref{prop3})
and the random error terms are, up to logarithmic factors, of the order $O(\frac{1}{\sqrt{n}})$
with respect to the sample size $n.$ 

The next result provides a more subtle oracle inequality that is 
akin to approximation error bounds of Proposition \ref{prop1}.
In this oracle inequality, the approximation error term due to von Neumann entropy penalization is $a^2(\log S)\eps^2$ (as in Proposition 
\ref{prop1}), so, it is of the order $O(\eps^2).$ Note that it is 
assumed implicitly that $a^2(\log S)<+\infty,$ i.e., that 
$S$ is of full rank and the matrix $\log S$ is well defined.
The random error terms are of the order $O(n^{-1})$ as $n\to\infty $ (up to logarithmic factors) with an exception of the term 
$\sigma_{\xi} (\sigma_X\vee \|{\mathbb E}X\|)
\|P_{L^{\perp}}SP_{L^{\perp}}\|_1
\sqrt{\frac{t_{m}}{n}},$
which depends on how well the oracle $S$ is approximated 
by low rank matrices. If $\|P_{L^{\perp}}SP_{L^{\perp}}\|_1$
is small, say of the order $n^{-1/2}$ for a subspace $L$
of a small dimension $r,$ this term becomes comparable 
to other terms in the bound, or even smaller. 
The inequalities hold only for the values 
of regularization parameter $\eps$ above certain threshold 
(so, this result does not apply to the simple least squares 
estimator). The first bound shows that if there is an oracle 
$S\in {\cal S}$ such that: (a) it is ``well aligned'', that is, 
$a(\log S)$ is small; (b) there exists a subspace $L$ of small 
dimension $r$ such that the oracle matrix $S$ is ``almost supported''
in $L,$ that is, $\|P_{L^{\perp}}SP_{L^{\perp}}\|_1$ is small; 
and (c) $S$ provides a good approximation of the density matrix 
$\rho,$ that is, $\|S-\rho\|_{L_2(\Pi)}^2$ is small, then 
the empirical solution $\hat \rho^{\eps}$ will be in the 
intersection of the $L_2(\Pi)$-ball and the Kullback-Leibler 
``ball'' of small enough radii around the oracle $S.$ 
The second bound is an oracle inequality showing how the 
$L_2(\Pi)$-error $\|\hat \rho^{\eps}-\rho\|_{L_2(\Pi)}^2$
depends on the properties of the oracle $S.$

\begin{theorem} 
\label{oracle}
There exist numerical constants $C>0, D>0$ such that the following 
holds. 
For all $t>0,$ for all $\lambda>0,$
for all 
$
\eps \geq 
D\eps_{n,m},
$
for all subspaces $L\subset {\mathbb C}^m$ with ${\rm dim}(L):=r,$
and for all $S\in {\cal S},$
with probability at least $1-e^{-t},$ 
\begin{eqnarray}
\label{bound_oracle_1}
&&
\|\hat \rho^{\eps}-S\|_{L_2(\Pi)}^2+\frac{\eps}{4}K(\hat \rho^{\eps};S)\leq 
(1+\lambda)\|S-\rho\|_{L_2(\Pi)}^2 + 
\frac{C}{\lambda}
\biggl[a^2(\log S)\eps^2\bigvee 
\\
&&
\nonumber
\sigma_{\xi}^2\beta^2(L) \frac{m r+\tau_{n}}{n} \bigvee 
\sigma_{\xi} (\sigma_X\vee \|{\mathbb E}X\|)
\|P_{L^{\perp}}SP_{L^{\perp}}\|_1
\sqrt{\frac{t_{m}}{n}}\bigvee
c_{\xi} U\frac{\tau_n\vee t_m}{n}\bigvee U^2 \frac{t_m}{n}\biggr]
\end{eqnarray}
and
\begin{eqnarray}
\label{bound_6'}
&&
\|\hat \rho^{\eps}-\rho\|_{L_2(\Pi)}^2
\leq (1+\lambda)\|S-\rho\|_{L_2(\Pi)}^2+
\frac{C}{\lambda}
\biggl[a^2(\log S)\eps^2\bigvee 
\sigma_{\xi}^2\beta^2(L) \frac{m r+\tau_n}{n} \bigvee 
\nonumber
\\
&&
\sigma_{\xi} (\sigma_X\vee \|{\mathbb E}X\|)
\|P_{L^{\perp}}SP_{L^{\perp}}\|_1
\sqrt{\frac{t_{m}}{n}}\bigvee
c_{\xi} U \frac{\tau_n\vee t_m}{n}\bigvee U^2\frac{t_m}{n}\biggr].
\end{eqnarray}
\end{theorem}


Next we give a version of (\ref{bound_oracle_1}) in a 
special case when $S=\rho^{\eps}.$ This provides bounds on 
random errors of estimation of the true penalized solution 
$\rho^{\eps}$ by its empirical version $\hat \rho^{\eps}$ both 
in the $L_2(\Pi)$ and in the Kullback-Leibler distances. 
Note that unlike the bounds for an arbitrary oracle $S,$
there is no dependence on the alignment coefficient $a(\log \rho^{\eps})$
in this case. The result essentially shows that as soon as the true solution $\rho^{\eps}$ is approximately low rank in the sense that $P_{L^{\perp}}\rho^{\eps}P_{L^{\perp}}$ is ``small'' for a subspace 
$L$ of a ``small'' dimension $r$ and $\rho^{\eps}$ provides a good 
approximation of the target density matrix $\rho,$ the empirical 
solution $\hat \rho^{\eps}$ would also provide a good approximation 
of $\rho$ and it would be approximately low rank.

\begin{theorem} 
\label{random_error_1}
There exist numerical constants $C>0, D>0$ such that the following 
holds. 
For all $t>0,$ 
for all 
$
\eps \geq 
D\eps_{n,m}
$
and for all subspaces $L\subset {\mathbb C}^m$ with ${\rm dim}(L):=r,$
with probability at least $1-e^{-t},$ 
\begin{eqnarray}
\label{bound_6}
&&
\|\hat \rho^{\eps}-\rho^{\eps}\|_{L_2(\Pi)}^2+ \eps K(\hat \rho^{\eps};\rho^{\eps})
\leq 
\nonumber
\\
&&
C
\biggl[\sigma_{\xi}^2\beta^2(L) \frac{m r+\tau_n}{n} \bigvee 
\sigma_{\xi} (\sigma_X\vee \|{\mathbb E}X\|)
\|P_{L^{\perp}}\rho^{\eps}P_{L^{\perp}}\|_1
\sqrt{\frac{t_{m}}{n}}\bigvee
\nonumber
\\
&&
U \|\rho^{\eps}-\rho\|_{L_2(\Pi)}
\sqrt{\frac{t_{m}}{n}}\bigvee
U^2\|\rho^{\eps}-\rho\|_1\frac{t_{m}}{n}\bigvee 
c_{\xi} U \frac{\tau_n\vee t_m}{n}\biggr].
\end{eqnarray}
\end{theorem}

{\bf Remark}. In the case when the noise is not necessarily bounded,
but $\|\xi\|_{\psi_1}<+\infty,$ the results still hold with the following 
simple modifications. In bounds (\ref{simple_1}), (\ref{simple_2}),
(\ref{simple_3}) and in the definition of $\eps_{n,m},$ the term $(c_{\xi}U\vee U^2)\frac{t_m}{n}$
is to be replaced by 
$$
\biggl(\|\xi\|_{\psi_1}U \log\biggl(\frac{\|\xi\|_{\psi_1}}{\sigma_{\xi}}\frac{U}{\sigma_X}\biggr)
\bigvee U^2\biggr)\frac{t_m}{n}.
$$
In the bounds of theorems \ref{oracle} and \ref{random_error_1}, 
the term $c_{\xi} U \frac{\tau_n\vee t_m}{n}$ is to be replaced 
by 
$$
\|\xi\|_{\psi_1}U \frac{\tau_n \log n}{n}\bigvee 
\|\xi\|_{\psi_1}U \log\biggl(\frac{\|\xi\|_{\psi_1}}{\sigma_{\xi}}\frac{U}{\sigma_X}\biggr)
\frac{t_m}{n}.
$$

We will provide a detailed proof of Theorem \ref{oracle}. The proof of Theorem \ref{simple} is its simplified version.
The proof of Theorem \ref{random_error_1}
relies on the bounds derived in the proof of Theorem \ref{oracle}. 
It is also possible to derive 
the oracle inequalities of Theorem \ref{oracle} 
from Theorem \ref{random_error_1} and from the approximation error bounds of Proposition \ref{prop1}. 
Throughout the proofs below, $C, C_1, \dots $ are numerical constants 
whose values might be different in different places. 

{\bf Proof of Theorem \ref{oracle}}.  
Denote 
$$
L_n(S):=n^{-1}\sum_{j=1}^n (Y_j-{\rm tr}(SX_j))^2 + 
\eps\ {\rm tr}(S\log S).
$$
For any $S\in {\cal S}$ of full rank and any direction $\nu\in {\mathbb M}_m({\mathbb C}),$ we have  
$$
DL_n(S;\nu) = 
2n^{-1}\sum_{j=1}^n(\langle S,X_j\rangle -Y_j)\langle \nu, X_j\rangle + \eps\ {\rm tr}(\nu \log S).
$$
By necessary 
conditions of extrema in the convex optimization problem (\ref{entropy_penalty}),
$
DL_n(\hat\rho^{\eps};\hat \rho^{\eps}-S)\leq 0,
$
which implies  
\begin{equation}
\label{DL-2}
DL(\hat \rho^{\eps};\hat \rho^{\eps}-S)-
DL(S;\hat \rho^{\eps}-S)
\leq  
-DL(S;\hat \rho^{\eps}-S)+
DL(\hat \rho^{\eps};\hat \rho^{\eps}-S)-
DL_n(\hat \rho^{\eps};\hat \rho^{\eps}-S).
\end{equation}
Note that 
$$
DL(\hat \rho^{\eps};\hat \rho^{\eps}-S)-
DL(S;\hat \rho^{\eps}-S)= 
2\|\hat \rho^{\eps}-S\|_{L_2(\Pi)}^2 + \eps K(\hat \rho^{\eps};S)
$$
(see the proof of Proposition \ref{prop3}) and 
$$DL(S;\hat \rho^{\eps}-S)= 
2\langle S-\rho,\hat \rho^{\eps}-S\rangle_{L_2(\Pi)}
+\eps {\rm tr}((\hat \rho^{\eps}-S)\log S).
$$
By a simple algebra similar to what has been already used in the 
proofs of propositions \ref{prop3}, \ref{prop1}, we get 
the following bound:
\begin{eqnarray}
\label{bound_1A}
&&
2\|\hat \rho^{\eps}-S\|_{L_2(\Pi)}^2+2\langle S-\rho,\hat \rho^{\eps}-S\rangle_{L_2(\Pi)}+
\eps K(\hat \rho^{\eps};S)=
\\
&&
\nonumber
\|\hat \rho^{\eps}-S\|_{L_2(\Pi)}^2+\|\hat \rho^{\eps}-\rho\|_{L_2(\Pi)}^2
-\|S-\rho\|_{L_2(\Pi)}^2+
\eps K(\hat \rho^{\eps};S)
\leq 
\\
&&
\nonumber
-\eps {\rm tr}((\hat \rho^{\eps}-S)\log S)
-\frac{2}{n} \sum_{j=1}^n \biggl(
\langle \hat \rho^{\eps}-S,X_j\rangle^2- 
{\mathbb E}\langle \hat \rho^{\eps}-S,X\rangle^2
\biggr)+ 
\\
&&
\nonumber
\frac{2}{n} \sum_{j=1}^n \biggl(
\langle S-\rho,X_j\rangle 
\langle \hat \rho^{\eps }-S,
X_j\rangle
- 
{\mathbb E}\langle S-\rho,X\rangle 
\langle \hat \rho^{\eps }-S,X\rangle\biggr)
-
\frac{2}{n}\sum_{j=1}^n \xi_j 
\langle \hat \rho^{\eps }-S,X_j\rangle.
\end{eqnarray}
Since 
$
\eps |{\rm tr}((\hat \rho^{\eps}-S)\log S)|\leq 
\eps a(\log S)\|\hat \rho^{\eps}-S\|_{L_2(\Pi)},
$ 
we get from (\ref{bound_1A}) that 
\begin{eqnarray}
\label{bound_1A_1}
&&
\|\hat \rho^{\eps}-S\|_{L_2(\Pi)}^2+\|\hat \rho^{\eps}-\rho\|_{L_2(\Pi)}^2
+\eps K(\hat \rho^{\eps};S)
\leq 
\\
&&
\nonumber
\|S-\rho\|_{L_2(\Pi)}^2 + 
\eps a(\log S)\|\hat \rho^{\eps}-S\|_{L_2(\Pi)}
-\frac{2}{n} \sum_{j=1}^n \biggl(
\langle \hat \rho^{\eps}-S,X_j\rangle^2- 
{\mathbb E}\langle \hat \rho^{\eps}-S,X\rangle^2
\biggr)+ 
\\
&&
\nonumber
\frac{2}{n} \sum_{j=1}^n \biggl(
\langle S-\rho,X_j\rangle 
\langle \hat \rho^{\eps }-S,
X_j\rangle
- 
{\mathbb E}\langle S-\rho,X\rangle 
\langle \hat \rho^{\eps }-S,X\rangle\biggr)
-
\frac{2}{n}\sum_{j=1}^n \xi_j 
\langle \hat \rho^{\eps }-S,X_j\rangle.
\end{eqnarray}

We need to bound the empirical processes in the right hand side of bound (\ref{bound_1A_1}). We will do it in several steps by bounding each term  separately.

{\bf Step 1}. To bound the first term note that
$$
\frac{1}{n} \sum_{j=1}^n \biggl(
\langle \hat \rho^{\eps}-S,X_j\rangle^2- 
{\mathbb E}\langle \hat \rho^{\eps}-S,X\rangle^2
\biggr) =
\biggl\langle 
(\hat \rho^{\eps}-S)\otimes (\hat \rho^{\eps}-S),
\frac{1}{n}\sum_{j=1}^n ((X_j\otimes X_j)-{\mathbb E}(X\otimes X))
\biggr\rangle .
$$
Therefore,
$$
\biggl|\frac{1}{n} \sum_{j=1}^n \biggl(
\langle \hat \rho^{\eps}-S,X_j\rangle^2- 
{\mathbb E}\langle \hat \rho^{\eps}-S,X\rangle^2
\biggr)\biggr| 
\leq \|\hat \rho^{\eps}-S\|_1^2
\biggl\|\frac{1}{n}\sum_{j=1}^n ((X_j\otimes X_j)-{\mathbb E}(X\otimes X))
\biggr\|.
$$
Note that $\|X\otimes X\|=\|X\|^2\leq U^2$ and also 
$
\|X\otimes X-{\mathbb E}(X\otimes X)\|\leq 2 U^2.
$
Using noncommutative Bernstein's inequality (see (\ref{bernstein-1})
in subsection 3.3) we can claim that with probability 
at least $1-e^{-t}$
$$
\biggl\|\frac{1}{n}\sum_{j=1}^n ((X_j\otimes X_j)-{\mathbb E}(X\otimes X))
\biggr\|
\leq 
4\biggl(\sigma_{X\otimes X}\sqrt{\frac{t+\log (2m^2)}{n}}\bigvee
U^2\frac{t+\log (2m^2)}{n} 
\biggr) 
$$
and, with the same probability,
$$
\biggl|\frac{1}{n} \sum_{j=1}^n \biggl(
\langle \hat \rho^{\eps}-S,X_j\rangle^2- 
{\mathbb E}\langle \hat \rho^{\eps}-S,X\rangle^2
\biggr)\biggr| 
\leq 
$$
$$
4\biggl(\sigma_{X\otimes X}\sqrt{\frac{t+\log (2m^2)}{n}}\bigvee
U^2\frac{t+\log (2m^2)}{n} 
\biggr) \|\hat \rho^{\eps}-S\|_1^2.
$$

{\bf Step 2}. The second term can be written as 
$$
\frac{1}{n} \sum_{j=1}^n \biggl(
\langle S-\rho,X_j\rangle 
\langle \hat \rho^{\eps }-S,
X_j\rangle
- 
{\mathbb E}\langle S-\rho,X\rangle 
\langle \hat \rho^{\eps }-S,X\rangle
\biggr)=
$$
$$
\biggl\langle \hat \rho^{\eps}-S,\frac{1}{n} \sum_{j=1}^n 
\biggl(
\langle S-\rho,X_j\rangle X_j
- 
{\mathbb E}\langle S-\rho,X\rangle X 
\biggr)\biggr\rangle 
$$
and bounded as follows
$$
\biggl|\frac{1}{n} \sum_{j=1}^n \biggl(
\langle S-\rho,X_j\rangle 
\langle \hat \rho^{\eps }-S,
X_j\rangle
- 
{\mathbb E}\langle S-\rho,X\rangle 
\langle \hat \rho^{\eps }-S,X\rangle
\biggr)\biggr|\leq 
$$
$$
\|\hat \rho^{\eps}-S\|_1\biggl\|\frac{1}{n} \sum_{j=1}^n 
\biggl(
\langle S-\rho,X_j\rangle X_j
- 
{\mathbb E}\langle S-\rho,X\rangle X 
\biggr)\biggr\|.
$$
We use again the noncommutative version of Bernstein's inequality 
to show that with probability at least $1-e^{-t}$
$$
\biggl\|\frac{1}{n} \sum_{j=1}^n 
\biggl(
\langle S-\rho,X_j\rangle X_j
- 
{\mathbb E}\langle S-\rho,X\rangle X 
\biggr)\biggr\|\leq 
$$
$$
4U \|S-\rho\|_{L_2(\Pi)}
\sqrt{\frac{t+\log (2m)}{n}}\bigvee
4 U^2\|S-\rho\|_1\frac{t+\log (2m)}{n},
$$
where we also used simple bounds 
$
\|{\mathbb E}
\langle S-\rho, X\rangle^2 X^2\|
\leq U^2 \|S-\rho\|_{L_2(\Pi)}^2
$
and 
$
\|\langle S-\rho, X\rangle X\|\leq U^2 \|S-\rho\|_1.
$
Since $\|\hat \rho^{\eps}-S\|_1\leq 2,$ we get 
$$
\biggl|\frac{1}{n} \sum_{j=1}^n \biggl(
\langle S-\rho,X_j\rangle 
\langle \hat \rho^{\eps }-S,
X_j\rangle
- 
{\mathbb E}\langle S-\rho,X\rangle 
\langle \hat \rho^{\eps }-S,X\rangle
\biggr)\biggr|\leq 
$$
$$
8U \|S-\rho\|_{L_2(\Pi)}
\sqrt{\frac{t+\log (2m)}{n}}\bigvee
8 U^2\|S-\rho\|_1\frac{t+\log (2m)}{n}.
$$

{\bf Step 3}. We turn now to bounding the third term in the right hand side of (\ref{bound_1A_1}). 
It is easy to decompose it as follows:
\begin{eqnarray}
\label{gugu}
&&
\frac{1}{n}\sum_{j=1}^n \xi_j 
\langle \hat \rho^{\eps }-S,X_j\rangle =
\biggl\langle P_{L^{\perp}}(\hat\rho^{\eps}-S)P_{L^{\perp}},
\frac{1}{n}\sum_{j=1}^n \xi_j P_{L^{\perp}}X_j P_{L^{\perp}}\biggr\rangle 
+
\nonumber
\\
&&
\frac{1}{n}\sum_{j=1}^n \xi_j 
\langle \hat \rho^{\eps }-S,{\cal P}_L X_j\rangle .
\end{eqnarray}
Note that 
$$
\biggl|\biggl\langle P_{L^{\perp}}(\hat\rho^{\eps}-S)P_{L^{\perp}},
\frac{1}{n}\sum_{j=1}^n \xi_j P_{L^{\perp}}X_jP_{L^{\perp}}\biggr\rangle
\biggr|
\leq \|P_{L^{\perp}}(\hat\rho^{\eps}-S)P_{L^{\perp}}\|_1
\biggl\|\frac{1}{n}\sum_{j=1}^n \xi_j P_{L^{\perp}}X_jP_{L^{\perp}}\biggr\|.
$$
Applying the noncommutative version of Bernstein's inequality one more time,
we have that with probability at least $1-e^{-t}$
$$  
\biggl\|\frac{1}{n}\sum_{j=1}^n \xi_j 
(P_{L^{\perp}}X_j P_{L^{\perp}}-
{\mathbb E}P_{L^{\perp}}XP_{L^{\perp}})\biggr\|\leq 
2\sigma_{\xi} \sigma_X 
\sqrt{\frac{t+\log (2m)}{n}}\bigvee
2 c_{\xi} U \frac{t+\log (2m)}{n},
$$
where we used a simple bound 
$
\|{\mathbb E}(P_{L^{\perp}} (X-{\mathbb E}X) P_{L^{\perp}})^2\|\leq \|{\mathbb E}(X-{\mathbb E}X)^2\|=\sigma_X^2.
$
Also, it follows from the classical Bernstein's inequality 
and the bound $\|{\mathbb E}(P_{L^{\perp}}XP_{L^{\perp}})\|\leq \|{\mathbb E}X\|$ that 
with probability at least $1-e^{-t}$
$$  
\biggl\|\frac{1}{n}\sum_{j=1}^n \xi_j 
{\mathbb E}P_{L^{\perp}}XP_{L^{\perp}}\biggr\|=
\biggl|\frac{1}{n}\sum_{j=1}^n \xi_j\biggr| 
\Bigl\|{\mathbb E}P_{L^{\perp}}XP_{L^{\perp}}\Bigr\|
\leq 
2\sigma_{\xi}\|{\mathbb E}X\| 
\sqrt{\frac{t}{n}}\bigvee
2 c_{\xi} \|\mathbb E X\| \frac{t}{n}.
$$
Hence, with probability at least $1-2e^{-t},$
$$
\biggl|\biggl\langle P_{L^{\perp}}(\hat\rho^{\eps}-S)P_{L^{\perp}},
\frac{1}{n}\sum_{j=1}^n \xi_j P_{L^{\perp}}X_j P_{L^{\perp}}\biggr\rangle
\biggr|\leq
$$
$$ 
2\|P_{L^{\perp}}(\hat\rho^{\eps}-S)P_{L^{\perp}}\|_1
\biggl[\sigma_{\xi} (\sigma_X+\|{\mathbb E}X\|)
\sqrt{\frac{t+\log (2m)}{n}}\bigvee
2c_{\xi} U \frac{t+\log (2m)}{n}\biggr].
$$
To bound the second term in the right hand side of (\ref{gugu}), denote 
$$
\alpha_n(\delta):=\sup_{\rho_1,\rho_2\in {\cal S}, \|\rho_1-\rho_2\|_{L_2(\Pi)}\leq \delta}\biggl|\frac{1}{n}\sum_{j=1}^n \xi_j 
\langle \rho_1-\rho_2,{\cal P}_L X_j\rangle\biggr| .
$$
Clearly, 
$
\biggl|\frac{1}{n}\sum_{j=1}^n \xi_j 
\langle \hat \rho^{\eps }-S,{\cal P}_L X_j\rangle \biggr|
\leq \alpha_n(\|\hat \rho^{\eps}-S\|_{L_2(\Pi)}).
$
To control $\alpha_n(\delta),$ we use Talagrand's 
concentration inequality for empirical processes.
It implies that, for all $\delta>0,$ with probability at least $1-e^{-s},$
\begin{equation}
\label{tal}
\alpha_n(\delta)\leq 2\biggl[{\mathbb E}\alpha_n(\delta)+
\sigma_{\xi} \beta (L)\delta \sqrt{\frac{s}{n}}+ 4 c_{\xi} U\frac{s}{n}
\biggr].
\end{equation}
Here we used the facts that 
$
{\mathbb E}\xi^2 \langle \rho_1-\rho_2, {\cal P}_L X\rangle^2
\leq \sigma_{\xi}^2 \beta^2(L) \|\rho_1-\rho_2\|_{L_2(\Pi)}^2  
$
and
$$
\Bigl|\xi \langle \rho_1-\rho_2, {\cal P}_L X\rangle\Bigr|\leq 
c_{\xi}\|\rho_1-\rho_2\|_1 \|{\cal P}_L X\|\leq 2c_{\xi}(\|X\|+\|P_{L^{\perp}} X P_{L^{\perp}}\|)
\leq 4c_{\xi}\|X\|\leq 4c_{\xi}U.
$$
We will make the bound on $\alpha_n(\delta)$ uniform 
in $\delta \in [Un^{-1}, 2U].$  
To this end, we apply bound (\ref{tal})
for $\delta=\delta_j=2^{-j+1}U,\ j=0,1,\dots $ and 
with $s=\tau_n:=t+ \log \log_2(2n) .$ The union bound and 
the monotonicity of $\alpha_n(\delta)$ with respect 
to $\delta$ implies that with probability at least 
$1-e^{-t}$ for all $\delta \in [Un^{-1},2U]$
\begin{equation}
\label{tal_1}
\alpha_n(\delta)\leq C\biggl[{\mathbb E}\alpha_n(\delta)+
\sigma_{\xi} \beta (L)\delta \sqrt{\frac{\tau_n}{n}}+ c_{\xi}U\frac{\tau_n}{n}
\biggr],
\end{equation}
where $C>0$ is a numerical constant.
Now it remains to bound the expected value ${\mathbb E}\alpha_n(\delta).$
Let $e_1,\dots, e_m$ be the orthonormal basis of ${\mathbb C}^m$ such 
that $L={\rm l.s.}\{e_1,\dots, e_r\}.$ Denote $E_{ij}(x)$ the entries 
of the linear transformation $x\in {\mathbb M}_m({\mathbb C})$ in 
this basis. Clearly, the function $\langle \rho_1-\rho_2, {\cal P}_L x\rangle$ belongs to the space ${\cal L}:={\rm l.s.}\{E_{ij}: i\leq r\ {\rm or}\ j\leq r\}$ of dimension $m^2-(m-r)^2=2 m r-r^2$ 
Therefore, 
$$
{\mathbb E}\alpha_n(\delta)\leq {\mathbb E}\sup_{f\in {\cal L},\|f\|_{L_2(\Pi)}\leq 
\beta(L)\delta}\biggl|\frac{2}{n}\sum_{j=1}^n \xi_j f(X_j)\biggr|.
$$
Using standard bounds for empirical processes indexed by 
finite dimensional function classes, we get 
$
{\mathbb E}\alpha_n(\delta)\leq
2\sqrt{2}
\sigma_{\xi}\beta(L)\delta \sqrt{\frac{m r}{n}}.
$
We can conclude that the following bound on $\alpha_n(\delta)$ holds 
with probability at least $1-e^{-t}$ for all $\delta\in [U n^{-1},2U]:$
\begin{equation}
\label{tal_2}
\alpha_n(\delta)\leq 
C\biggl[\sigma_{\xi}\beta(L)\delta \sqrt{\frac{m r}{n}}+
\sigma_{\xi} \beta (L)\delta \sqrt{\frac{\tau_n}{n}}+ c_{\xi}U\frac{\tau_n}{n}
\biggr].
\end{equation}
Note that since $\|\hat \rho^{\eps}-S\|_1\leq 2$ and $\|X\|\leq U,$
we have
$
\|\hat \rho^{\eps}-S\|_{L_2(\Pi)}^2 = 
{\mathbb E}\langle \hat \rho^{\eps}-S,X\rangle \leq 4U^2,
$
so, $\|\hat \rho^{\eps}-S\|_{L_2(\Pi)}\leq 2U.$
As a result, with probability at least $1-e^{-t},$
we either have $\|\hat \rho^{\eps}-S\|_{L_2(\Pi)}<Un^{-1},$ or 
$$
\biggl|\frac{1}{n}\sum_{j=1}^n \xi_j 
\langle \hat \rho^{\eps }-S,{\cal P}_L X_j\rangle \biggr|
\leq 
C\biggl[\sigma_{\xi}\beta(L)
\|\hat \rho^{\eps}-S\|_{L_2(\Pi)} 
\sqrt{\frac{m r}{n}}+
\sigma_{\xi} \beta (L) 
\|\hat \rho^{\eps}-S\|_{L_2(\Pi)} 
\sqrt{\frac{\tau_n}{n}}+ c_{\xi}U\frac{\tau_n}{n}
\biggr].
$$
In the first case, we still have 
$$
\biggl|\frac{1}{n}\sum_{j=1}^n \xi_j 
\langle \hat \rho^{\eps }-S,{\cal P}_L X_j\rangle \biggr|
\leq 
C\biggl[\sigma_{\xi}\beta(L)
\frac{U}{n} 
\sqrt{\frac{m r}{n}}+
\sigma_{\xi} \beta (L) 
\frac{U}{n} 
\sqrt{\frac{\tau_n}{n}}+ c_{\xi}U\frac{\tau_n}{n}
\biggr].
$$
Let us assume in what follows that 
$\|\hat \rho^{\eps}-S\|_{L_2(\Pi)}\geq Un^{-1}$
since another case is even easier to handle. 

We now substitute the bounds of steps 1--3 in the right hand side 
of (\ref{bound_1A_1}) to get the following inequality that holds 
with some constant $C>0$ and with probability 
at least $1-5e^{-t}:$
\begin{eqnarray}
\label{bound_2A}
&&
\|\hat \rho^{\eps}-S\|_{L_2(\Pi)}^2+ \|\hat \rho^{\eps}-\rho\|_{L_2(\Pi)}^2+ \eps K(\hat \rho^{\eps};S)
\leq 
\\
&&
\nonumber
\|S-\rho\|_{L_2(\Pi)}^2+
\eps a(\log S)\|\hat \rho^{\eps}-S\|_{L_2(\Pi)}+
\\
&&
\nonumber
16\biggl(\sigma_{X\otimes X}\sqrt{\frac{t_{m}}{n}}\bigvee
U^2\frac{t_{m}}{n} 
\biggr) \|\hat \rho^{\eps}-S\|_1^2+
16U \|S-\rho\|_{L_2(\Pi)}
\sqrt{\frac{t_{m}}{n}}\bigvee 16U^2\frac{t_m}{n}+
\\
&&
\nonumber
4\|P_{L^{\perp}}(\hat\rho^{\eps}-S)P_{L^{\perp}}\|_1
\biggl[\sigma_{\xi} (\sigma_X+\|{\mathbb E}X\|) 
\sqrt{\frac{t_{m}}{n}}\bigvee
2c_{\xi} U \frac{t_{m}}{n}\biggr]+
\\
&&
\nonumber 
C\biggl[\sigma_{\xi}\beta(L)
\|\hat \rho^{\eps}-S\|_{L_2(\Pi)} 
\sqrt{\frac{m r+\tau_n}{n}}
\bigvee 
c_{\xi}U\frac{\tau_n}{n}
\biggr].
\end{eqnarray}
Under the assumption 
$\eps \geq D\eps_{n,m}$ with a sufficiently large 
constant $D>0,$ 
it is easy to get that 
\begin{equation}
\label{zeta1}
16\biggl(\sigma_{X\otimes X}\sqrt{\frac{t_m}{n}}\bigvee
U^2\frac{t_m}{n}\biggr) \|\hat \rho^{\eps}-S\|_1^2
\leq 
\frac{\eps}{2}\|\hat \rho^{\eps}-S\|_1^2 \leq 
\frac{\eps}{2}K(\hat \rho^{\eps};S).
\end{equation}
Also, by Proposition \ref{rank_00}, 
$$
\|P_{L^{\perp}}(\hat\rho^{\eps}-S)P_{L^{\perp}}\|_1
\leq 
\|P_{L^{\perp}}\hat \rho^{\eps}P_{L^{\perp}}\|_1+
\|P_{L^{\perp}} S P_{L^{\perp}}\|_1\leq 
3\|P_{L^{\perp}}S P_{L^{\perp}}\|_1+2K(\hat \rho^{\eps};S),
$$
and, under the same assumption 
that $\eps \geq D\eps_{n,m}$ with a sufficiently large 
constant $D>0,$ 
\begin{eqnarray}
\label{zeta2}
&&
4
\|P_{L^{\perp}}(\hat\rho^{\eps}-S)P_{L^{\perp}}\|_1
\biggl[\sigma_{\xi} (\sigma_X+\|{\mathbb E}X\|)
\sqrt{\frac{t_m}{n}}\bigvee
2c_{\xi} U \frac{t_m}{n}\biggr]
\leq 
\\
&&
C\|P_{L^{\perp}} S P_{L^{\perp}}\|_1
\biggl[\sigma_{\xi} (\sigma_X\vee \|{\mathbb E}X\|)
\sqrt{\frac{t_m}{n}}\bigvee
c_{\xi} U \frac{t_m}{n}\biggr]
+\frac{\eps}{4} K(\hat \rho^{\eps};S).
\nonumber
\end{eqnarray}
Combining bounds (\ref{zeta1}) and (\ref{zeta2}) with (\ref{bound_2A}) yields 
\begin{eqnarray}
\label{bound_2AB}
&&
\|\hat \rho^{\eps}-S\|_{L_2(\Pi)}^2+ \|\hat \rho^{\eps}-\rho\|_{L_2(\Pi)}^2
+\frac{\eps}{4} K(\hat \rho^{\eps};S)
\leq 
\\
&&
\nonumber
\|S-\rho\|_{L_2(\Pi)}^2+
\eps a(\log S)\|\hat \rho^{\eps}-S\|_{L_2(\Pi)}+
\\
&&
\nonumber
C\biggl[
\|\hat \rho^{\eps}-S\|_{L_2(\Pi)} 
\sigma_{\xi}\beta(L)
\sqrt{\frac{m r+\tau_n}{n}}\bigvee 
U \|S-\rho\|_{L_2(\Pi)}
\sqrt{\frac{t_m}{n}}\bigvee 
\\
&&
\nonumber
\|P_{L^{\perp}} S P_{L^{\perp}}\|_1
\sigma_{\xi} (\sigma_X\vee \|{\mathbb E}X\|)
\sqrt{\frac{t_m}{n}}
\bigvee 
c_{\xi}U\frac{\tau_n\vee t_m}{n}\bigvee U^2\frac{t_m}{n}\biggr]
\end{eqnarray}
with some constant $C>0.$ It follows from the last inequality that 
\begin{equation}
\label{dada}
\|\hat \rho^{\eps}-S\|_{L_2(\Pi)}^2\leq 
A\|\hat \rho^{\eps}-S\|_{L_2(\Pi)}+B-\frac{\eps}{4}K(\hat \rho^{\eps};S),
\end{equation}
where 
$
A:=\frac{\eps}{2} a(\log S)+ C\sigma_{\xi}\beta(L)
\sqrt{\frac{m r+\tau_n}{n}}
$
and 
\begin{eqnarray}
&&
\nonumber
B:=\|S-\rho\|_{L_2(\Pi)}^2 - \|\hat \rho^{\eps}-\rho\|_{L_2(\Pi)}^2+ 
\\
&&
\nonumber
C\biggl[\|S-\rho\|_{L_2(\Pi)}
U\sqrt{\frac{t_m}{n}}\bigvee \|P_{L^{\perp}} S P_{L^{\perp}}\|_1
\sigma_{\xi} (\sigma_X\vee \|{\mathbb E}X\|)
\sqrt{\frac{t_m}{n}}
\bigvee c_{\xi}U\frac{\tau_n\vee t_m}{n}\bigvee U^2\frac{t_m}{n}
\biggr].
\end{eqnarray}
It is easy to check that 
$$
\|\hat \rho^{\eps}-S\|_{L_2(\Pi)}^2
\leq \biggl(\frac{A+\sqrt{A^2+4(B-(\eps/4)K(\hat \rho^{\eps};S))}}{2}\biggr)^2\leq
\biggl(A+\sqrt{\biggl(B-\frac{\eps}{4}K(\hat \rho^{\eps};S)\biggr)_{+}}\biggr)^2.
$$
If
$
\frac{\eps}{4}K(\hat \rho^{\eps};S)\geq B,
$ 
then $\|\hat \rho^{\eps}-S\|_{L_2(\Pi)}^2\leq A^2,$
which, in view of (\ref{dada}), implies 
$$
\|\hat \rho^{\eps}-S\|_{L_2(\Pi)}^2+\frac{\eps}{4}K(\hat \rho^{\eps};S)
\leq A^2+B.
$$
Otherwise, we have 
$
\|\hat \rho^{\eps}-S\|_{L_2(\Pi)}^2
\leq 
A^2+2A\sqrt{B}+ B-\frac{\eps}{4}K(\hat \rho^{\eps};S),
$
which, for all $\lambda>0,$ implies 
$$
\|\hat \rho^{\eps}-S\|_{L_2(\Pi)}^2+
\frac{\eps}{4}K(\hat \rho^{\eps};S)
\leq (\frac{2}{\lambda}+1)A^2 + (1+\lambda/2)B.
$$
In both cases, by the definitions of $A$ and $B$ and by elementary algebra, one can easily get the bound 
\begin{eqnarray}
&&
\|\hat \rho^{\eps}-S\|_{L_2(\Pi)}^2+
\|\hat \rho^{\eps}-\rho\|_{L_2(\Pi)}^2+
\frac{\eps}{4}K(\hat \rho^{\eps};S)\leq 
\nonumber
\\
&&
(1+\lambda)\|S-\rho\|_{L_2(\Pi)}^2 + 
\frac{C}{\lambda}
\biggl[a^2(\log S)\eps^2\bigvee 
\sigma_{\xi}^2\beta^2(L) \frac{m r+\tau_n}{n} \bigvee 
\nonumber
\\
&&
\sigma_{\xi} (\sigma_X\vee \|{\mathbb E}X\|)
\|P_{L^{\perp}}SP_{L^{\perp}}\|_1
\sqrt{\frac{t_{m}}{n}}\bigvee
c_{\xi} U \frac{\tau_n\vee t_m}{n}\bigvee U^2\frac{t_m}{n}\biggr]
\end{eqnarray}
that holds with probability at least $1-5e^{-t}$
and with a sufficiently large constant $C.$
To replace the probability $1-5e^{-t}$ by $1-e^{-t},$
it is enough to replace $t$ by $t+\log 5$ and to adjust 
the values of constants $C, D$ accordingly.

\qed

{\bf Proof of Theorem \ref{simple}}. 
We get back to bound (\ref{bound_1A}) in the proof 
of Theorem \ref{oracle}. This time, we bound the term 
${\rm tr}((\hat \rho^{\eps}-S)\log S)$ 
in (\ref{bound_1A}) in a slightly different way
$$
|{\rm tr}((\hat \rho^{\eps}-S)\log S)|
\leq \|\log S\|\|\hat \rho- S\|_1 \leq 2\|\log S\|,
$$
which leads to the following bound (instead of bound (\ref{bound_1A_1})):
\begin{eqnarray}
\label{bound_1A_BB}
&&
\|\hat \rho^{\eps}-S\|_{L_2(\Pi)}^2+ \|\hat \rho^{\eps}-\rho\|_{L_2(\Pi)}^2
+
\eps K(\hat \rho^{\eps};S)
\leq 
\|S-\rho\|_{L_2(\Pi)}^2+
\eps \|\log S)\|+
\\
&&
\nonumber
-
\frac{1}{n} \sum_{j=1}^n \biggl(
\langle \hat \rho^{\eps}-S,X_j\rangle^2- 
{\mathbb E}\langle \hat \rho^{\eps}-S,X\rangle^2
\biggr)+ 
\\
&&
\nonumber
\frac{1}{n} \sum_{j=1}^n \biggl(
\langle S-\rho,X_j\rangle 
\langle \hat \rho^{\eps }-S,
X_j\rangle
- 
{\mathbb E}\langle S-\rho,X\rangle 
\langle \hat \rho^{\eps }-S,X\rangle\biggr)
-
\frac{1}{n}\sum_{j=1}^n \xi_j 
\langle \hat \rho^{\eps }-S,X_j\rangle.
\end{eqnarray}
To bound the empirical processes in the right hand side, we again use 
the bounds of steps 1--3 in the proof of Theorem \ref{oracle}. 
The bound of Step 1 yields 
$$
\biggl|\frac{1}{n} \sum_{j=1}^n \biggl(
\langle \hat \rho^{\eps}-S,X_j\rangle^2- 
{\mathbb E}\langle \hat \rho^{\eps}-S,X\rangle^2
\biggr)\biggr| 
\leq 
16\biggl(\sigma_{X\otimes X}\sqrt{\frac{t+\log (2m^2)}{n}}\bigvee
U^2\frac{t+\log (2m^2)}{n} 
\biggr)
$$
and it follows from the bound of Step 2 that 
$$
\biggl|\frac{1}{n} \sum_{j=1}^n \biggl(
\langle S-\rho,X_j\rangle 
\langle \hat \rho^{\eps }-S,
X_j\rangle
- 
{\mathbb E}\langle S-\rho,X\rangle 
\langle \hat \rho^{\eps }-S,X\rangle
\biggr)\biggr|\leq 
$$
$$
8U \|S-\rho\|_{L_2(\Pi)}
\sqrt{\frac{t+\log (2m)}{n}}\bigvee
16 U^2\frac{t+\log (2m)}{n}.
$$
Instead of more complicated derivation of Step 3, we now use 
noncommutative and classical Bernstein's inequalities to get
that with probability at least $1-2 e^{-t}$ 
$$
\biggl|n^{-1}\sum_{j=1}^n \xi_j \langle \hat \rho^{\eps}-S,X_j\rangle\biggr|
\leq \|\hat \rho^{\eps}-S\|_1 
\biggl\|n^{-1}\sum_{j=1}^n \xi_j X_j\biggr\|\leq 
2
\biggl\|n^{-1}\sum_{j=1}^n \xi_j (X_j-{\mathbb E}X)\biggr\|
+
$$
$$
2\|{\mathbb E}X\| \biggl|n^{-1}\sum_{j=1}^n \xi_j \biggr|\leq 
4(\sigma_{\xi}\sigma_X+\|{\mathbb E}X\|)\sqrt{\frac{t+\log(2m)}{n}}
\bigvee 12c_{\xi}U\frac{t+\log(2m)}{n}.
$$
Using these inequalities, we derive from (\ref{bound_1A_BB}) that 
with some numerical constant $C>0$ and with probability at least 
$1-4e^{-t},$
\begin{eqnarray}
\label{bound_1A_BBB}
&&
\|\hat \rho^{\eps}-S\|_{L_2(\Pi)}^2+ 
\|\hat \rho^{\eps}-\rho\|_{L_2(\Pi)}^2
+\eps K(\hat \rho^{\eps};S)
\leq 
\|S-\rho\|_{L_2(\Pi)}^2+
\eps \|\log S\|+
\\
&&
\nonumber
+ 
C\biggl[
U \|S-\rho\|_{L_2(\Pi)}\sqrt{\frac{t_m}{n}}+
(\sigma_{X\otimes X}\vee \sigma_{\xi}\sigma_X\vee \|{\mathbb E}X\|)\sqrt{\frac{t_m}{n}}\bigvee (c_{\xi}U\vee U^2)\frac{t_m}{n}
\biggr],
\end{eqnarray}
which implies the result in the case when $\|\log S\|\leq \log \Gamma.$
To finish the proof, it is enough, given an arbitrary $S\in {\cal S}$ 
(even such that $\log S$ does not exist), to apply bound (\ref{bound_1A_BBB}) 
to $S_{\delta}=(1-\delta)S+ \delta \frac{I_m}{m},$ where $\delta\in (0,1).$
Clearly, $\|\log S_{\delta}\|\leq \log\frac{m}{\delta}$ and we also have 
$
\|S-S_{\delta}\|_{L_2(\Pi)}^2\leq 4 \delta^2 {\mathbb E}\|X\|^2
$
(see the proof of Proposition \ref{low_rank_approx}). Taking 
$
\delta:= \frac{\sqrt{\eps}}{{\mathbb E}^{1/2}\|X\|^2}\wedge 1,
$
it is easy to complete the proof in the case when $\|\log S\|\geq \log \Gamma .$

\qed

{\bf Proof of Theorem \ref{random_error_1}}. Note that similarly
to $\rho^{\eps},$ $\hat \rho^{\eps}$ is also a matrix of full 
rank and $\log \hat \rho^{\eps}$ is well defined. 
By necessary 
conditions of extrema in convex problems (\ref{entropy_penalty}) and 
(\ref{entropy_penalty_true}),
we have  
$
DL_n(\hat\rho^{\eps};\hat \rho^{\eps}-\rho^{\eps})\leq 0
$
and  
$
DL(\rho^{\eps};\hat \rho^{\eps}-\rho^{\eps})\geq 0.
$
Subtracting the second inequality from the first one yields 
\begin{equation}
\label{DL-1}
DL(\hat \rho^{\eps};\hat \rho^{\eps}-\rho^{\eps})-
DL(\rho^{\eps};\hat \rho^{\eps}-\rho^{\eps})
\leq 
DL(\hat \rho^{\eps};\hat \rho^{\eps}-\rho^{\eps})-
DL_n(\hat \rho^{\eps};\hat \rho^{\eps}-\rho^{\eps}).
\end{equation}
By a simple algebra already used in the proof of Theorem \ref{oracle}, 
this easily leads to the following bound:
$$
2\|\hat \rho^{\eps}-\rho^{\eps}\|_{L_2(\Pi)}^2+ \eps K(\hat \rho^{\eps};\rho^{\eps})
\leq 2{\mathbb E}\langle \hat \rho^{\eps}-\rho, X\rangle 
\langle \hat \rho^{\eps}-\rho^{\eps},X\rangle
-
2n^{-1}\sum_{j=1}^n(\langle \hat \rho^{\eps},X_j\rangle -Y_j)
\langle \hat \rho^{\eps}-\rho^{\eps}, X_j\rangle,
$$
which can be further rewritten as 
\begin{eqnarray}
\label{bound_1}
&&
2\|\hat \rho^{\eps}-\rho^{\eps}\|_{L_2(\Pi)}^2+ \eps K(\hat \rho^{\eps};\rho^{\eps})
\leq 
-
\frac{2}{n} \sum_{j=1}^n \biggl(
\langle \hat \rho^{\eps}-\rho^{\eps},X_j\rangle^2- 
{\mathbb E}\langle \hat \rho^{\eps}-\rho^{\eps},X\rangle^2
\biggr) -
\\
&&
\nonumber
\frac{2}{n} \sum_{j=1}^n \biggl(
\langle \rho^{\eps}-\rho,X_j\rangle 
\langle \hat \rho^{\eps }-\rho^{\eps},
X_j\rangle
- 
{\mathbb E}\langle \rho^{\eps}-\rho,X\rangle 
\langle \hat \rho^{\eps }-\rho^{\eps},X\rangle
\biggr)
-\frac{2}{n}\sum_{j=1}^n \xi_j 
\langle \hat \rho^{\eps }-\rho^{\eps},X_j\rangle .
\end{eqnarray}

We use the bounds of steps 1--3 of the proof of Theorem \ref{oracle}  with $S=\rho^{\eps}$ to control each term in the right hand side of (\ref{bound_1}). Substituting these bounds in  
(\ref{bound_1}), we get the following inequality
that holds with probability at least $1-5e^{-t}:$
\begin{eqnarray}
\label{bound_2}
&&
2\|\hat \rho^{\eps}-\rho^{\eps}\|_{L_2(\Pi)}^2+ \eps K(\hat \rho^{\eps};\rho^{\eps})
\leq 
\\
&&
\nonumber
8\biggl(\sigma_{X\otimes X}\sqrt{\frac{t+\log (2m^2)}{n}}\bigvee
U^2\frac{t+\log (2m^2)}{n} 
\biggr) \|\hat \rho^{\eps}-\rho^{\eps}\|_1^2+
\\
&&
\nonumber
16U \|\rho^{\eps}-\rho\|_{L_2(\Pi)}
\sqrt{\frac{t+\log (2m)}{n}}\bigvee
16U^2\|\rho^{\eps}-\rho\|_1\frac{t+\log (2m)}{n}+
\\
&&
\nonumber
4\|P_{L^{\perp}}(\hat\rho^{\eps}-\rho^{\eps})P_{L^{\perp}}\|_1
\biggl[\sigma_{\xi} (\sigma_X+\|{\mathbb E}X\|)
\sqrt{\frac{t+\log (2m)}{n}}\bigvee
2c_{\xi} U \frac{t+\log (2m)}{n}\biggr]+
\\
&&
\nonumber 
C\biggl[\sigma_{\xi}\beta(L)
\|\hat \rho^{\eps}-\rho^{\eps}\|_{L_2(\Pi)} 
\sqrt{\frac{m r}{n}}+
\sigma_{\xi} \beta (L) 
\|\hat \rho^{\eps}-\rho^{\eps}\|_{L_2(\Pi)} 
\sqrt{\frac{\tau_n}{n}}+ c_{\xi}U\frac{\tau_n}{n}
\biggr].
\end{eqnarray}
Arguing exactly as in the proof of Theorem \ref{oracle}, we can simplify 
(\ref{bound_2}) to get 
\begin{eqnarray}
\label{bound_4}
&&
2\|\hat \rho^{\eps}-\rho^{\eps}\|_{L_2(\Pi)}^2+ \frac{\eps}{4} K(\hat \rho^{\eps};\rho^{\eps})
\leq 
\\
&&
\nonumber
16U \|\rho^{\eps}-\rho\|_{L_2(\Pi)}
\sqrt{\frac{t+\log (2m)}{n}}\bigvee
16 U^2\|\rho^{\eps}-\rho\|_1\frac{t+\log (2m)}{n}+
\\
&&
\nonumber
12
\|P_{L^{\perp}}\rho^{\eps}P_{L^{\perp}}\|_1
\biggl[\sigma_{\xi} (\sigma_X+\|{\mathbb E}X\|)
\sqrt{\frac{t+\log (2m)}{n}}\bigvee
2c_{\xi} U \frac{t+\log (2m)}{n}\biggr]+
\\
&&
\nonumber 
C\biggl[\sigma_{\xi}\beta(L)
\|\hat \rho^{\eps}-\rho^{\eps}\|_{L_2(\Pi)} 
\sqrt{\frac{m r}{n}}+
\sigma_{\xi} \beta (L) 
\|\hat \rho^{\eps}-\rho^{\eps}\|_{L_2(\Pi)} 
\sqrt{\frac{\tau_n}{n}}+ c_{\xi}U\frac{\tau_n}{n}
\biggr].
\end{eqnarray}
It is easy now to solve this for $\|\hat \rho^{\eps}-\rho^{\eps}\|_{L_2(\Pi)}$
and to derive the following explicit bound 
on the random error that holds with probability 
at least $1-5e^{-t}$ and with some numerical 
constant $C>0:$
\begin{eqnarray}
\label{bound_5}
&&
\|\hat \rho^{\eps}-\rho^{\eps}\|_{L_2(\Pi)}^2+ \eps K(\hat \rho^{\eps};\rho^{\eps})
\leq 
C
\biggl[\sigma_{\xi}^2\beta^2(L) \frac{m r+\tau_n}{n} 
\bigvee c_{\xi}U\frac{\tau_n}{n} \bigvee 
\\
&&
\nonumber
U \|\rho^{\eps}-\rho\|_{L_2(\Pi)}
\sqrt{\frac{t+\log (2m)}{n}}\bigvee
U^2\|\rho^{\eps}-\rho\|_1\frac{t+\log (2m)}{n}\bigvee 
\\
&&
\nonumber
\|P_{L^{\perp}}\rho^{\eps}P_{L^{\perp}}\|_1
\biggl(\sigma_{\xi} (\sigma_X\vee \|{\mathbb E}X\|) 
\sqrt{\frac{t+\log (2m)}{n}}\bigvee
c_{\xi} U \frac{t+\log (2m)}{n}\biggr)\biggr],
\end{eqnarray}
which easily implies the result.

\qed

{\bf Example 1. Matrix completion (continuation)}. Recall that, in this example, 
$\{e_i: i=1,\dots, m\}$ is the canonical basis of ${\mathbb C}^m$ 
and the following set of 
Hermitian matrices forms an orthonormal basis of ${\mathbb M}_m({\mathbb C})$ (the matrix completion basis):
$$
\Bigl\{e_i \otimes e_i: i=1,\dots, m\Bigr\}
\bigcup 
\biggl\{\frac{1}{\sqrt{2}}(e_i\otimes e_j+e_j\otimes e_i): 1\leq i<j\leq m\biggr\}
$$
$$
\bigcup 
\biggl\{\frac{i}{\sqrt{2}}(e_i\otimes e_j-e_j\otimes e_i): 1\leq i<j\leq m\biggr\}.
$$  
Assume that $X$ is sampled at random from this basis. Recall that in 
this case, for all matrices $A,$
$
\|A\|_{L_2(\Pi)}^2 = m^{-2}\|A\|_2^2.
$
Obviously, 
$
\|e_i\otimes e_i\|=1,\ i=1,\dots, m
$
and, for all $i<j,$
$$
\biggl\|\frac{1}{\sqrt{2}}(e_i\otimes e_j+e_j\otimes e_i)\biggr\|=\frac{1}{\sqrt{2}}, \ \ 
\biggl\|\frac{i}{\sqrt{2}}(e_i\otimes e_j-e_j\otimes e_i)\biggr\|=
\frac{1}{\sqrt{2}}.
$$
Therefore, $\|X\|\leq U=1.$
We also have 
$$
\sigma_X^2 \leq  \|{\mathbb E}X^2\|= 
\sup_{v\in {\mathbb C}^m,|v|=1}
{\mathbb E}\langle X^2 v,v\rangle=
\sup_{v\in {\mathbb C}^m,|v|=1}
{\mathbb E}\langle Xv, Xv\rangle=
\sup_{v\in {\mathbb C}^m,|v|=1}
{\mathbb E}|Xv|^2.
$$
Note that, if $X=e_i\otimes e_i, i=1,\dots, m,$ then 
$ 
|Xv|^2 = |e_i\langle e_i,v\rangle|^2=|\langle e_i,v\rangle|^2. 
$
If $X=\frac{1}{\sqrt{2}}(e_i\otimes e_j+e_j\otimes e_i), i<j,$ then  
$$
|Xv|^2= \frac{1}{2}
|e_i \langle e_j,v\rangle + e_j \langle e_i,v\rangle|^2= 
\frac{1}{2} 
\Bigl(|\langle e_j,v\rangle|^2 + |\langle e_i,v\rangle|^2\Bigr)
$$
and, similarly, if $X=\frac{i}{\sqrt{2}}(e_i\otimes e_j-e_j\otimes e_i),
i<j,$
then also 
$
|Xv|^2=
\frac{1}{2} \Bigl(|\langle e_j,v\rangle|^2 + |\langle e_i,v\rangle|^2\Bigr).
$
Therefore, for $|v|=1,$
$$
{\mathbb E}|Xv|^2 = m^{-2}\sum_{i=1}^m |\langle e_i,v\rangle|^2
+2 m^{-2} \frac{1}{2}\sum_{i<j}\Bigl(|\langle e_j,v\rangle|^2 + |\langle e_i,v\rangle|^2\Bigr) \leq  
$$
$$
m^{-2}|v|^2 + m^{-2}m (|v|^2+|v|^2)\leq 3m^{-1},
$$
which implies that 
$
\sigma_X\leq \frac{\sqrt{3}}{\sqrt{m}}.
$
By a similar simple computation,
$
\sigma_{X\otimes X}\leq \frac{4}{\sqrt{m}} .
$
Now we can derive the following corollary of Theorem \ref{oracle}.
Let 
$$
\eps_{n,m}:=
(\sigma_{\xi}m^{-1/2}\vee m^{-1/2})\sqrt{\frac{t_m}{n}}\bigvee
(c_{\xi}\vee 1)\frac{t_m}{n}
$$
and let $\eps =D\eps_{n,m}$ for a sufficiently large constant $D>0.$

\begin{corollary} 
\label{oracle-1}
There exists a numerical constant $C>0$ such that the following 
holds. For all $t>0,$ for all $\lambda>0,$ for all sufficiently 
large $D$ and for $\eps=D\eps_{n,m},$
for all matrices $S\in {\cal S}$ of rank $r,$ 
with probability at least $1-e^{-t},$ 
\begin{eqnarray}
\label{bound_6'''}
&&
\|\hat \rho^{\eps}-\rho\|_{L_2(\Pi)}^2
\leq (1+\lambda) \|S-\rho\|_{L_2(\Pi)}^2+
\frac{C}{\lambda}
\biggl[
D^2\biggl((\sigma_{\xi}^2 \vee 1)\frac{r m t_m}{n}\bigvee
\nonumber
\\
&&
(c_{\xi}^2\vee 1)\frac{r m^2 t_m^2}{n^2}\biggr)\log^2 (m n)
\bigvee 
\sigma_{\xi}^2 \frac{\tau_n}{n} \bigvee 
c_{\xi}\frac{\tau_n\vee t_m}{n}\vee \frac{t_m}{n}\biggr].
\end{eqnarray}
\end{corollary}

{\bf Proof}. First observe that for all matrices $S\in {\cal S}$
of full rank (for which $\log S$ exists) and for all subspaces 
$L\subset {\mathbb C}^m$ with ${\rm dim}(L)=r,$ we have, with 
probability at least $1-e^{-t}$ and with an arbitrary $\lambda>0$ 
\begin{eqnarray}
\label{bound_6''}
&&
\|\hat \rho^{\eps}-\rho\|_{L_2(\Pi)}^2
\leq (1+\lambda/2)\|S-\rho\|_{L_2(\Pi)}^2+
\frac{2C}{\lambda}
\biggl[a^2(\log S)\biggl((\sigma_{\xi}^2\vee 1)\frac{t_m}{mn}+(c_{\xi}^2\vee 1)
\frac{t_m^2}{n}\biggr)  
\bigvee 
\nonumber
\\
&&
\sigma_{\xi}^2 \frac{m r+\tau_n}{n} \bigvee 
\sigma_{\xi} m^{-1/2}
\|P_{L^{\perp}}SP_{L^{\perp}}\|_1
\sqrt{\frac{t_m}{n}}\bigvee
c_{\xi}\frac{\tau_n\vee t_m}{n}\bigvee \frac{t_m}{n}\biggr].
\end{eqnarray}
This immediately follows from Theorem \ref{oracle} since,
in the case under consideration,  
$\beta(L)=1,$
$\sigma_X \leq 3^{1/2} m^{-1/2},$ $\sigma_{X\otimes X}\leq 4 m^{-1/2},$
$U=1.$ Note also that in this case $\Lambda(L)=m$ (recall the definition
of $\Lambda (L)$ given before Proposition \ref{low_rank_approx})
and 
$$
a(\log S)\leq m\inf_{c}\|\log S+cI_m\|_2.
$$
Suppose now that $S\in {\cal S}$ is an arbitrary oracle of rank $r.$
Then there exists a subspace $L$ of dimension $r$ such that 
$P_{L^{\perp}}SP_{L^{\perp}}=0.$ We will use bound (\ref{bound_6''})
for 
$S_{\delta}:=(1-\delta)S+\delta \frac{I_m}{m},$
where 
$
\delta =\eps \wedge 1,
$ 
as we did in the 
proof of Proposition 
\ref{low_rank_approx}.
As in this proof, we have, for some constant $C_1>0,$ 
$$
a(\log S_{\delta})
\leq 
m\sqrt{r}\log \biggl(1+\frac{m}{\delta}\biggr)
\leq C_1 m\sqrt{r}\log (mn)
$$
and 
$$
\|S-S_{\delta}\|_{L_2(\Pi)}^2 
\leq 4 \delta^2 {\mathbb E}\|X\|^2 \leq 4 \delta^2\leq 4\eps^2.
$$
Finally, note that 
$$
\|P_{L^{\perp}}S_{\delta}P_{L^{\perp}}\|_1 \leq 
(1-\delta)\|P_{L^{\perp}}S P_{L^{\perp}}\|_1+
\delta \|P_{L^{\perp}}(I_m/m)P_{L^{\perp}}\|_1\leq 
\delta \leq \eps .
$$
Substituting these bounds in (\ref{bound_6''}) (with $S$ replaced by $S_{\delta}$) and bounding 
$\|S_{\delta}-\rho\|_{L_2(\Pi)}^2$ in terms of $\|S-\rho\|_{L_2(\Pi)}^2$ and $\|S_{\delta}-S\|_{L_2(\Pi)}^2$
(similarly to what was done in the proof of Proposition 
\ref{low_rank_approx}), it is easy to derive (\ref{bound_6'''}) from
(\ref{bound_6''}). 
Note that we can drop the term 
$\sigma_{\xi}^2 \frac{m r}{n}$ since it is dominated 
by $(\sigma_{\xi}^2 \vee 1)\frac{r m t_{m}}{n}\log^2(mn).$

\qed

Similarly, it is easy to obtain another corollary where the $L_2(\Pi)$-error of estimator $\hat \rho^{\eps}$ is controlled 
in terms of Gibbs oracles. Recall the notations at the end of 
Section 4 and also denote 
$
\Gamma_r:= \|H_{\leq r}\|_2^2 = \sum_{k=1}^r \gamma_k^2.
$ 

\begin{corollary} 
\label{oracle-Gibbs}
There exists a numerical constant $C>0$ such that the following 
holds. For all $t>0,$ for all $\lambda>0,$ for all sufficiently 
large $D$ and for $\eps=D\eps_{n,m},$
for all Hermitian matrices $H$ and for all $r\leq m,$ 
with probability at least $1-e^{-t},$ 
\begin{eqnarray}
\label{bound_6-Gibbs}
&&
\|\hat \rho^{\eps}-\rho\|_{L_2(\Pi)}^2
\leq (1+\lambda) \|\rho_{H}-\rho\|_{L_2(\Pi)}^2+
\frac{C}{\lambda}
\biggl[\frac{\delta_r^2(H)}{m^2}\bigvee 
D^2\biggl((\sigma_{\xi}^2 \vee 1)\frac{\Gamma_r m t_m}{n}\bigvee
\nonumber
\\
&&
(c_{\xi}^2\vee 1)\frac{\Gamma_r m^2 t_m^2}{n^2}\biggr)
\bigvee 
\sigma_{\xi}^2 \frac{m r+\tau_n}{n} \bigvee 
c_{\xi} \frac{\tau_n\vee t_m}{n}\bigvee \frac{t_m}{n}
\biggr].
\end{eqnarray}
\end{corollary}

{\bf Example 2. Pauli basis (continuation)}. We now turn to another example described 
in the Introduction, the example of the Pauli basis. Recall that in this 
case $m=2^k$ and we are considering the basis of the space ${\mathbb M}_{2^k}({\mathbb C})$ that consists of all matrices of the form 
$W_{i_1}\otimes \dots \otimes W_{i_k},$ $W_i=\frac{1}{\sqrt{2}}\sigma_i, i=1,\dots, 4$ being normalized $2\times 2$ Pauli matrices.
Note that $\|W_i\|_2=1$ and $\|W_i\|=\frac{1}{\sqrt{2}}.$ 
The design variable $X$ is picked at random from this basis. 
We still have 
$
\|A\|_{L_2(\Pi)}^2 = m^{-2}\|A\|_2^2.
$
However, now  
$$
\|W_{i_1}\otimes \dots \otimes W_{i_k}\|= \|W_{i_1}\|\dots \|W_{i_k}\|
=\biggl(\frac{1}{\sqrt{2}}\biggr)^k=2^{-k/2}=m^{-1/2} 
$$
implying that $\|X\|=m^{-1/2}$ and $U=m^{-1/2}.$ 
To state a corollary of Theorem \ref{oracle} in this case,
we take $\eps:=D\eps_{n,m},$ where 
$$
\eps_{n,m}:=
(\sigma_{\xi}m^{-1/2}\vee m^{-1})\sqrt{\frac{t_m}{n}}\bigvee
(c_{\xi}m^{-1/2}\vee m^{-1})\frac{t_m}{n}.
$$

The following results are similar to corollaries \ref{oracle-1}
and \ref{oracle-Gibbs}.

\begin{corollary} 
\label{oracle-2}
There exists a numerical constant $C>0$ such that the following 
holds. 
For all $t>0,$ for all $\lambda>0,$ for all sufficiently large 
$D>0$ and for $\eps=D\eps_{n,m},$
for all matrices $S\in {\cal S}$ of rank $r,$ 
with probability at least $1-e^{-t},$
\begin{eqnarray}
\label{bound_6'''A}
&&
\|\hat \rho^{\eps}-\rho\|_{L_2(\Pi)}^2
\leq (1+\lambda)\|S-\rho\|_{L_2(\Pi)}^2+
\frac{C}{\lambda}
\biggl[D^2
\biggl((\sigma_{\xi}^2 \vee m^{-1})\frac{r m t_m}{n}\bigvee
\nonumber
\\
&&
(c_{\xi}^2\vee m^{-1})\frac{r m t_m^2}{n^2}\biggr)\log^2 (m n)
\bigvee 
\sigma_{\xi}^2 \frac{\tau_n}{n} 
\bigvee
c_{\xi}m^{-1/2}\frac{\tau_n\vee t_m}{n}\bigvee \frac{t_m}{m n}\biggr].
\end{eqnarray}
\end{corollary}

\begin{corollary} 
\label{oracle-Gibbs-2}
There exists a numerical constant $C>0$ such that the following 
holds. For all $t>0,$ for all $\lambda>0,$ for all sufficiently 
large $D$ and for $\eps=D\eps_{n,m},$
for all Hermitian matrices $H$ and for all $r\leq m,$ 
with probability at least $1-e^{-t},$ 
\begin{eqnarray}
\label{bound_6-Gibbs-2}
&&
\|\hat \rho^{\eps}-\rho\|_{L_2(\Pi)}^2
\leq (1+\lambda) \|\rho_{H}-\rho\|_{L_2(\Pi)}^2+
\frac{C}{\lambda}
\biggl[\frac{\delta_r^2(H)}{m^2}\bigvee 
D^2\biggl((\sigma_{\xi}^2 \vee m^{-1})\frac{\Gamma_r m t_m}{n}\bigvee
\nonumber
\\
&&
(c_{\xi}^2\vee m^{-1})\frac{\Gamma_r m^2 t_m^2}{n^2}\biggr)
\bigvee 
\sigma_{\xi}^2 \frac{m r+\tau_n}{n} \bigvee 
c_{\xi}m^{-1/2}\frac{\tau_n\vee t_m}{n}
\bigvee \frac{t_m}{m n}
\biggr].
\end{eqnarray}
\end{corollary}

Note that the bounds of corollaries
\ref{oracle-1}-\ref{oracle-Gibbs-2}  
can be also proved in the case when the noise is unbounded, in particular, 
Gaussian (see the remark after Theorem \ref{random_error_1}). For the Pauli basis, 
this immediately leads to Theorem \ref{intro_th_1}
stated in the Introduction.

\section{Oracle Inequalities: Subgaussian Design Case}

In this section, we turn to the case of \it subgaussian design matrices. \rm More precisely, we assume that $X$ is a Hermitian random matrix 
with distribution $\Pi$ such that, for some constant $b_0>0$ and for all 
Hermitian matrices $A\in {\mathbb M}_m({\mathbb C}),$ 
$\langle A, X\rangle$ is a subgaussian random variable with 
parameter $b_0\|A\|_{L_2(\Pi)}.$
This implies that ${\mathbb E}X=0$ and, for some constant $b_1>0,$ 
\begin{equation}
\label{SG2}
\Bigl\|\langle A,X\rangle\Bigr\|_{\psi_2}\leq b_1 \|A\|_{L_2(\Pi)},\ 
A\in {\mathbb M}_m({\mathbb C}).
\end{equation}
In addition to this, assume that, for some constant $b_2>0,$
\begin{equation}
\label{SG2-1}
\|A\|_{L_2(\Pi)}=\Bigl\|\langle A,X\rangle\Bigr\|_{L_2(\Pi)}\leq b_2 \|A\|_2,\ 
A\in {\mathbb M}_m({\mathbb C}).
\end{equation}
A Hermitian random matrix $X$ satisfying the above conditions will be called a \it subgaussian \rm matrix.
Moreover, if $X$ also satisfies the condition 
\begin{equation}
\label{SG1}
\|A\|_{L_2(\Pi)}^2={\mathbb E}|\langle A,X\rangle|^2=\|A\|_2^2,\ A\in {\mathbb M}_m({\mathbb C}),
\end{equation}
then it will be called an \it isotropic subgaussian \rm matrix.  
As it was already mentioned in the introduction, the last class of matrices includes such examples as Gaussian and Rademacher design matrices. 
It easily follows from the basic properties of Orlicz norms 
(see, e.g., van der Vaart and Wellner (1996), p. 95) that 
for subgaussian matrices 
$
\|A\|_{L_p(\Pi)}={\mathbb E}^{1/p}
\Bigl|\langle A,X\rangle\Bigr|^p\leq c_p b_1 b_2 \|A\|_2^2
$
and 
$
\|A\|_{\psi_1}:=\Bigl\|\langle A,X\rangle\Bigr\|_{\psi_1}
\leq c b_1 b_2 \|A\|_2,
A\in {\mathbb M}_m({\mathbb C}), p\geq 1,
$
with some numerical constants $c_p>0$ and $c>0.$

The following is a version of a well known fact (see, e.g., Rudelson and Vershynin (2010), Proposition 2.4).

\begin{proposition}
\label{subgauss}
Let $X$ be a subgaussian $m\times m$ matrix. Then, there exists a 
constant $B>0$ such that 
$$
\Bigl\|\|X\|\Bigr\|_{\psi_2}\leq B\sqrt{m}.
$$ 
\end{proposition}

{\bf Proof}. Let $M\subset S^{m-1}:=\{u\in {\mathbb C}^m: |u|=1\}$ be an  $\eps$-net of the unit sphere in ${\mathbb C}^m$ of the smallest cardinality.
It is easy to see that ${\rm card}(M)\leq (1+2/\eps)^m$ and  
$$
\|X\|=
\sup_{u,v\in S^{m-1}}\langle X u,v\rangle
\leq (1-\eps)^{-2}\max_{u,v\in M}\langle X u,v\rangle .
$$
Take $\eps=1/2.$ Using standard bounds for Orlicz norms of a maximum
(see, e.g., van der Vaart and Wellner (1996), Lemma 2.2.2), we get that,
with some constants $C_1,C_2,B>0,$
$$
\Bigl\|\|X\|\Bigr\|_{\psi_2} \leq 
4 \Bigl\|\max_{u,v\in M}\langle X u,v\rangle \Bigr\|_{\psi_2}
\leq 
C_1 \psi_2^{-1}({\rm card}^2(M)) 
\max_{u,v\in M}
\Bigl\|\langle X u,v\rangle\Bigr\|_{\psi_2}\leq 
$$
$$
C_2 \sqrt{\log {\rm card}(M)} 
\max_{u,v\in M}\Bigl\|\langle X, u\otimes v\rangle\Bigr\|_{\psi_2}\leq 
C_2 \sqrt{\log {\rm card}(M)} 
\max_{u,v\in M}\|u\otimes v\|_2\leq B\sqrt{m}.
$$
\qed 

Below, we give oracle inequalities and random error bounds  
in the subgaussian design case. 
We will use the following notations. Given $t>0,$
let 
$$
t_m:=t+\log (2m), \ \ \tau_n:=t+\log \log_2(2n),\ \ {\rm and}\ \ 
t_{n,m}:= \tau_n\log n \vee t_m. 
$$
Also, denote
$
c_{\xi}:= \|\xi\|_{\psi_2}\log \frac{\|\xi\|_{\psi_2}}{\sigma_{\xi}}
$
and let 
$$
\eps_{n,m}:=\sigma_{\xi} 
\sqrt{\frac{mt_{m}}{n}}\bigvee
c_{\xi}\frac{\sqrt{m}t_{m}}{n}
$$
(clearly, we assume here that the noise has a bounded $\psi_2$-norm).

\begin{theorem}
\label{simple-SG}
There exist constants $C>0, c>0$ such that the following 
holds. For all $t>0$ and $\lambda>0$ 
such that $\tau_n\leq c \lambda^2 n,$
for all $S\in {\cal S}$
and for all $\eps\in [0,1],$
with probability at least $1-e^{-t}$
\begin{eqnarray}
\label{simple_1-SG}
&&
\nonumber
\|\hat \rho^{\eps}-S\|_{L_2(\Pi)}^2\leq (1+\lambda)
\|S-\rho\|_{L_2(\Pi)}^2 
+C\biggl[\eps 
\biggl(\|\log S\|\bigwedge 
\log \frac{m}{\eps}\biggr)
\bigvee 
\sigma_{\xi}\sqrt{\frac{mt_m}{n}}\bigvee 
\\
&&
\frac{m t_m}{n\lambda } \bigvee 
(c_{\xi}\vee \sqrt{m})\frac{\sqrt{m} t_{n,m}}{n} 
\biggr]
\end{eqnarray}
and 
\begin{eqnarray}
\label{simple_2-SG}
&&
\nonumber
\|\hat \rho^{\eps}-\rho\|_{L_2(\Pi)}^2\leq 
(1+\lambda)\|S-\rho\|_{L_2(\Pi)}^2 
+C\biggl[\eps 
\biggl(\|\log S\|\wedge 
\log \frac{m}{\eps}\biggr)
\bigvee 
\sigma_{\xi}\sqrt{\frac{m t_m}{n}}\bigvee
\\
&& 
\frac{m t_m}{n\lambda } \bigvee 
(c_{\xi}\vee \sqrt{m})\frac{\sqrt{m} t_{n,m}}{n} 
\biggr].
\end{eqnarray}
In particular,
\begin{equation}
\label{simple_3-SG}
\nonumber
\|\hat \rho^{\eps}-\rho\|_{L_2(\Pi)}^2\leq 
C\biggl[\eps 
\biggl(\|\log \rho\|\wedge 
\log \frac{m}{\eps}\biggr)\bigvee 
\sigma_{\xi}\sqrt{\frac{m t_m}{n}}\bigvee 
(c_{\xi}\vee \sqrt{m})\frac{\sqrt{m} t_{n,m}}{n} 
\biggr].
\end{equation}
\end{theorem}


We now turn to more subtle oracle inequalities that take into 
account low rank properties of oracles $S\in {\cal S}.$

\begin{theorem} 
\label{oracle-SG}
There exist numerical constants $C>0, D>0, c>0$ such that the following 
holds. 
For all $t>0$ and $\lambda>0$ such that $\tau_n\leq c \lambda^2 n,$
for all 
$
\eps \geq 
D\eps_{n,m},
$
for all subspaces $L\subset {\mathbb C}^m$ with ${\rm dim}(L):=r$ 
and for all $S\in {\cal S},$
with probability at least $1-e^{-t},$ 
\begin{eqnarray}
\label{bound_oracle_1SG}
&&
\|\hat \rho^{\eps}-S\|_{L_2(\Pi)}^2+\frac{\eps}{4}K(\hat \rho^{\eps};S)\leq 
(1+\lambda)\|S-\rho\|_{L_2(\Pi)}^2 + 
\\
&&
\frac{C}{\lambda}
\biggl[a^2(\log S)\eps^2\bigvee 
\sigma_{\xi}^2 \beta^2 (L)\frac{m r+\tau_n}{n} \bigvee 
\sigma_{\xi}
\|P_{L^{\perp}}SP_{L^{\perp}}\|_1
\sqrt{\frac{m t_m}{n}}\bigvee
(c_{\xi} \vee \sqrt{m}) \frac{\sqrt{m} t_{n,m}}{n}\biggr]
\nonumber
\end{eqnarray}
and
\begin{eqnarray}
\label{bound_6'SG}
&&
\|\hat \rho^{\eps}-\rho\|_{L_2(\Pi)}^2
\leq (1+\lambda)\|S-\rho\|_{L_2(\Pi)}^2+
\frac{C}{\lambda}
\biggl[a^2(\log S)\eps^2\bigvee 
\sigma_{\xi}^2 \frac{m r+\tau_n}{n} \bigvee 
\nonumber
\\
&&
\sigma_{\xi} 
\|P_{L^{\perp}}SP_{L^{\perp}}\|_1
\sqrt{\frac{m t_{m}}{n}}\bigvee
(c_{\xi} \vee \sqrt{m}) \frac{\sqrt{m} t_{n,m}}{n}\biggr].
\end{eqnarray}
\end{theorem}

Similarly to the previous section, we also derived bounds on 
the random error $\|\hat \rho^{\eps}-\rho^{\eps}\|_{L_2(\Pi)}^2.$

\begin{theorem} 
\label{random_error_SG}
There exist numerical constants $C>0, D>0, c>0$ such that the following 
holds. Under the assumption that $\tau_n\leq cn,$
for all $t>0,$ for all 
$
\eps \geq 
D\eps_{n,m}
$
and for all subspaces $L\subset {\mathbb C}^m$ with ${\rm dim}(L):=r,$
with probability at least $1-e^{-t},$ 
\begin{eqnarray}
\label{bound_6_SG}
&&
\|\hat \rho^{\eps}-\rho^{\eps}\|_{L_2(\Pi)}^2+ \eps K(\hat \rho^{\eps};\rho^{\eps})
\leq 
C
\biggl[\sigma_{\xi}^2 \beta^2(L)\frac{m r+\tau_n}{n} \bigvee 
\sigma_{\xi} 
\|P_{L^{\perp}}\rho^{\eps}P_{L^{\perp}}\|_1
\sqrt{\frac{m t_m}{n}}\bigvee
\nonumber
\\
&&
\|\rho^{\eps}-\rho\|_{L_2(\Pi)}
\sqrt{\frac{m t_m}{n}}\bigvee 
(c_{\xi}\vee \sqrt{m})  \frac{\sqrt{m} t_{n,m}}{n}\biggr].
\end{eqnarray}
\end{theorem}

We will give only the proof of Theorem \ref{oracle-SG}.

{\bf Proof}. 
It follows the lines of the proof of Theorem \ref{oracle} very closely.
The main changes are in the bounds of steps 1--3 of this proof that have 
to be modified in the subgaussian design case. The rest of the proof is 
straightforward. 

In Step 1, we have to bound the following quantity: 
$$
\frac{1}{n} \sum_{j=1}^n \biggl(
\langle \hat \rho^{\eps}-S,X_j\rangle^2- 
{\mathbb E}\langle \hat \rho^{\eps}-S,X\rangle^2
\biggr).
$$
To this end, we will study the empirical process 
$$
\Delta_{n}(\delta):=\sup_{f\in {\cal F}_{\delta}}\biggl|
n^{-1}\sum_{j=1}^n (f^2 (X_j)-Pf^2)
\biggr|,
$$
where 
$
{\cal F}_{\delta}:= \{\langle S_1-S_2,\cdot\rangle :S_1,S_2
\in {\cal S}, \|S_1-S_2\|_{L_2(\Pi)}\leq \delta\}.
$
Clearly, 
$$
\biggl|\frac{1}{n} \sum_{j=1}^n \biggl(
\langle \hat \rho^{\eps}-S,X_j\rangle^2- 
{\mathbb E}\langle \hat \rho^{\eps}-S,X\rangle^2
\biggr)\biggr|
\leq \Delta_n(\|\hat \rho^{\eps}-S\|_{L_2(\Pi)}).
$$
Our goal is to obtain an upper bound on $\Delta_n(\delta)$ 
uniformly in $\delta \in [(m/n)^{1/2}, 2b_2].$ First we use a 
version of Talagrand's concentration inequality for empirical 
processes indexed by unbounded functions due to Adamczak (see subsection 3.2).
It implies that with some constant $C>0$ and with probability at least $1-e^{-t}$ 
\begin{equation}
\label{adam}
\Delta_n(\delta)\leq 2 {\mathbb E}\Delta_n(\delta)+ C\delta^2 \sqrt{\frac{t}{n}} + C\frac{m t\log n}{n}.
\end{equation}
Here we used the following bounds on the uniform variance and on the 
envelope of the function class ${\cal F}_{\delta}^2:$ for the uniform 
variance, with some constant $c>0,$
$$
\sup_{f\in {\cal F}_{\delta}}(Pf^4)^{1/2}= 
\sup_{S_1,S_2\in {\cal S},\|S_1-S_2\|_{L_2(\Pi)}\leq 
\delta}
{\mathbb E}^{1/2}\langle S_1-S_2,X\rangle^4=
$$
$$
\sup_{S_1,S_2\in {\cal S},\|S_1-S_2\|_{L_2(\Pi)}\leq 
\delta}\|S_1-S_2\|_{L_4(\Pi)}^2 \leq c \delta^2,
$$
by the equivalence properties of the norms in Orlicz spaces.
For the envelope, 
$$
\sup_{f\in {\cal F}_{\delta}}f^2(X)=
\sup_{S_1,S_2\in {\cal S},\|S_1-S_2\|_{L_2(\Pi)}\leq 
\delta}
\langle S_1-S_2,X\rangle^2\leq 4 \|X\|^2 
$$
and 
$$
\Bigl\|\max_{1\leq i\leq n}\sup_{f\in {\cal F}_{\delta}}f^2(X_i)\Bigr\|_{\psi_1}\leq 
c_1\Bigl\|\|X\|^2\Bigr\|_{\psi_1}\log n 
\leq  c_2 \Bigl\|\|X\|\Bigr\|_{\psi_2}^2 \log n
\leq c_3 m \log n,
$$
for some constants $c_1,c_2,c_3>0,$ where we used well known 
inequalities for maxima of random variables in Orlicz spaces 
(see, e.g., Lemma 2.2.2 in van der Vaart and Wellner (1996)). 

To bound the expectation ${\mathbb E} \Delta_n(\delta)$ we use a 
recent result by Mendelson (2010) 
(see subsection 3.2; in fact, even earlier 
result by Klartag and Mendelson (2005) with the $\psi_2$-diameter 
instead of $\psi_1$-diameter would suffice for our 
purposes). It gives 
\begin{equation}
\label{shahar}
{\mathbb E}\Delta_n(\delta)\leq c \biggl[
\sup_{f\in {\cal F}_{\delta}}\|f\|_{\psi_1}
\frac{\gamma_2({\cal F}_{\delta};\psi_2)}{\sqrt{n}}
\bigvee \frac{\gamma_2^2({\cal F}_{\delta};\psi_2)}{n}
\biggr]
\end{equation}
with some constant $c>0.$
It follows from (\ref{SG2})
that the $\psi_1$ and $\psi_2$-norms  
of functions from the class ${\cal F}_{\delta}$ can be bounded from above by a constant 
times the $L_2(P)$-norm. As a result,
\begin{equation}
\label{diam}
\sup_{f\in {\cal F}_{\delta}}\|f\|_{\psi_1}\leq c\delta 
\end{equation}
and the following bound holds for 
Talagrand's generic chaining complexities:
\begin{equation}
\label{gamma_2}
\gamma_2({\cal F}_{\delta};\psi_2)\leq 
\gamma_2({\cal F}_{\delta};c\|\cdot\|_{L_2(\Pi)}),
\end{equation}
where $c$ is a constant. Let $G$ be a symmetric real valued random 
matrix with independent centered Gaussian entries $\{g_{ij}\}$ on 
the diagonal and above, where ${\mathbb E}g_{ii}^2=1$ and 
${\mathbb E}g_{ij}^2=\frac{1}{2}, i\neq j.$ Then, using condition  (\ref{SG2-1}),
we have that, 
for some constant $c_1>0,$ 
$$
{\mathbb E}|\langle S_1, G\rangle - \langle S_2, G\rangle|^2 =
\|S_1-S_2\|_2^2\geq c_1\|S_1-S_2\|_{L_2(\Pi)}^2,
$$ 
and it easily follows from Talagrand's generic chaining 
bound that, for some constant $C>0,$ 
\begin{equation}
\label{tala}
\gamma_2({\cal F}_{\delta};c\|\cdot\|_{L_2(\Pi)})
\leq C 
{\mathbb E}
\sup_{S_1,S_2\in {\cal S},\|S_1-S_2\|_{L_2(\Pi)}\leq 
\delta}|\langle S_1-S_2, G\rangle|=: C \omega (G;\delta).
\end{equation}
It follows from (\ref{shahar}), (\ref{diam}), (\ref{gamma_2}) and 
(\ref{tala}) that 
\begin{equation}
\label{shahar_1}
{\mathbb E}\Delta_n(\delta)\leq C 
\biggl[
\delta
\frac{\omega(G;\delta)}{\sqrt{n}}
\bigvee
\frac{\omega^2(G;\delta)}{n}
\biggr].
\end{equation}
To bound 
$
{\mathbb E}
\sup_{S_1,S_2\in {\cal S},\|S_1-S_2\|_{L_2(\Pi)}\leq 
\delta}
|\langle S_1-S_2, G\rangle|,
$
note that 
$$
\Bigl|\langle S_1-S_2, G\rangle\Bigr|\leq \|S_1-S_2\|_1 \|G\|\leq 2\|G\|, 
$$
and, by Proposition \ref{subgauss},
$$
\omega (G;\delta)={\mathbb E}
\sup_{\rho_1,\rho_2\in {\cal S},\|\rho_1-\rho_2\|_{L_2(\Pi)}\leq 
\delta}\Bigl|\langle S_1-S_2, G\rangle 
\Bigr|\leq 2{\mathbb E}\|G\|\leq c\sqrt{m}.
$$
Substituting this bound in (\ref{shahar_1}) yields that, 
for some constant $C>0,$
\begin{equation}
\label{shahar_2}
{\mathbb E}\Delta_n(\delta)\leq C \biggl[
\delta\sqrt{\frac{m}{n}}
\bigvee \frac{m}{n}
\biggr]
\end{equation}
and combining (\ref{shahar_2}) with (\ref{adam}) gives that 
with probability at least $1-e^{-t}$ 
\begin{equation}
\label{adam1}
\Delta_n(\delta)\leq C\biggl[\delta\sqrt{\frac{m}{n}}
\bigvee \frac{m}{n}\bigvee 
\delta^2\sqrt{\frac{t}{n}} \bigvee \frac{m t\log n}{n}\biggr].
\end{equation}
It is easy to make bound (\ref{adam1}) uniform in 
$\delta\in [(m/n)^{1/2},2b_2]$
by a simple discretization argument (as we did in Step 3 of the proof 
of Theorem \ref{oracle}). This leads to the following result:
with probability at least $1-e^{-t},$ for all $\delta\in [(m/n)^{1/2},2b_2],$
\begin{equation}
\label{adam2}
\Delta_n(\delta)\leq 
C\biggl[\delta\sqrt{\frac{m}{n}}
\bigvee \frac{m}{n}\bigvee
\delta^2\sqrt{\frac{\tau_n}{n}} \bigvee 
\frac{m \tau_n \log n}{n}\biggr],
\end{equation}
where $\tau_n=t+\log\log_2(2n).$
Thus, with the same probability and with a proper choice of constant 
$C>0$
$$
\biggl|\frac{1}{n} \sum_{j=1}^n \biggl(
\langle \hat \rho^{\eps}-S,X_j\rangle^2- 
{\mathbb E}\langle \hat \rho^{\eps}-S,X\rangle^2
\biggr)\biggr|\leq 
$$
$$
C\biggl[\|\hat \rho^{\eps}-S\|_{L_2(\Pi)}
\sqrt{\frac{m}{n}}
\bigvee \frac{m}{n}\bigvee 
\|\hat \rho^{\eps}-S\|_{L_2(\Pi)}^2\sqrt{\frac{\tau_n}{n}} \bigvee 
\frac{m\tau_n\log n}{n}\biggr]
$$
provided that $\|\hat \rho^{\eps}-S\|_{L_2(\Pi)}\in [(m/n)^{1/2}, 2b_2].$

Similarly to Step 2 of the proof of Theorem \ref{oracle}, 
we have to bound the expression  
$$
\frac{1}{n} \sum_{j=1}^n \biggl(
\langle S-\rho,X_j\rangle 
\langle \hat \rho^{\eps }-S,
X_j\rangle
- 
{\mathbb E}\langle S-\rho,X\rangle 
\langle \hat \rho^{\eps }-S,X\rangle
\biggr)=
$$
$$
\biggl\langle \hat \rho^{\eps}-S,\frac{1}{n} \sum_{j=1}^n 
\biggl(
\langle S-\rho,X_j\rangle X_j
- 
{\mathbb E}\langle S-\rho,X\rangle X 
\biggr)\biggr\rangle .
$$
We use the bound 
$$
\biggl|\frac{1}{n} \sum_{j=1}^n \biggl(
\langle S-\rho,X_j\rangle 
\langle \hat \rho^{\eps }-S,
X_j\rangle
- 
{\mathbb E}\langle S-\rho,X\rangle 
\langle \hat \rho^{\eps }-S,X\rangle
\biggr)\biggr|\leq 
$$
$$
\|\hat \rho^{\eps}-S\|_1
\biggl\|\frac{1}{n} \sum_{j=1}^n 
\biggl(
\langle S-\rho,X_j\rangle X_j
- 
{\mathbb E}\langle S-\rho,X\rangle X 
\biggr)\biggr\|.
$$
and Proposition \ref{bernstein_A} with $\alpha=1.$
Note that 
$$
\|{\mathbb E}
\langle S-\rho, X\rangle^2 X^2\|
\leq 
{\mathbb E}\langle S-\rho, X\rangle^2 \|X\|^2\leq 
{\mathbb E}^{1/2}\langle S-\rho, X\rangle^4 {\mathbb E}^{1/2}\|X\|^4 \leq c m\|S-\rho\|_{L_2(\Pi)}^2
$$
with a constant $c>0.$ Also,
$$
\Bigl\|\|\langle S-\rho, X\rangle X\|\|_{\psi_1}=  
\Bigl\||\langle S-\rho, X\rangle| \|X\|\Bigr\|_{\psi_1}\leq
c_1 \|\langle S-\rho, X\rangle\|_{\psi_2}\Bigl\|\|X\|\Bigr\|_{\psi_2}
\leq c_2 \sqrt{m}\|S-\rho\|_{L_2(\Pi)}
$$
with some constants $c_1,c_2>0.$
Finally, note that 
$$
\|\hat \rho^{\eps}-S\|_{L_2(\Pi)}\leq b_2 \|S-\rho\|_2 
\leq b_2 \|S-\rho\|_1 \|S-\rho\|\leq 4b_2,
$$ 
since, for $S,\rho\in {\cal S},$ $\|S-\rho\|_1\leq 2$ and $\|S-\rho\|\leq 2.$  
Using the fact $\|\hat \rho^{\eps}-S\|_1\leq 2,$
Proposition \ref{bernstein_A} and the previous bounds imply that with probability at least $1-e^{-t}$ and with some constants $C_1, C_2, C>0,$
$$
\biggl|\frac{1}{n} \sum_{j=1}^n \biggl(
\langle S-\rho,X_j\rangle 
\langle \hat \rho^{\eps }-S,
X_j\rangle
- 
{\mathbb E}\langle S-\rho,X\rangle 
\langle \hat \rho^{\eps }-S,X\rangle
\biggr)\biggr|\leq 
$$
$$\|\hat \rho^{\eps}-S\|_1 \biggl\|\frac{1}{n} \sum_{j=1}^n 
\biggl(
\langle S-\rho,X_j\rangle X_j
- 
{\mathbb E}\langle S-\rho,X\rangle X 
\biggr)\biggr\|\leq 
$$
$$
C_1\biggl[\|S-\rho\|_{L_2(\Pi)}
\sqrt{\frac{m(t+\log (2m))}{n}}\bigvee
\frac{\sqrt{m}(t+\log (2m))}{n}
\|S-\rho\|_{L_2(\Pi)}\log\frac{C_1 \sqrt{m}\|S-\rho\|_{L_2(\Pi)}}{\sqrt{m}\|S-\rho\|_{L_2(\Pi)}}\biggr]\leq 
$$
$$
C\biggl[\|S-\rho\|_{L_2(\Pi)}
\sqrt{\frac{m(t+\log (2m))}{n}}\bigvee
\frac{\sqrt{m}(t+\log (2m))}{n}
\biggr].
$$

We now modify the bounds of Step 3 of the proof of Theorem \ref{oracle}. We need to bound the following 
expression:
$$
\frac{1}{n}\sum_{j=1}^n \xi_j 
\langle \hat \rho^{\eps }-S,X_j\rangle =
\biggl\langle P_{L^{\perp}}(\hat\rho^{\eps}-S)P_{L^{\perp}},
\frac{1}{n}\sum_{j=1}^n \xi_j P_{L^{\perp}}X_j P_{L^{\perp}}\biggr\rangle 
+
\frac{1}{n}\sum_{j=1}^n \xi_j 
\langle \hat \rho^{\eps }-S,{\cal P}_L X_j\rangle .
$$
As in the proof of Theorem \ref{oracle}, 
$$
\biggl|\biggl\langle P_{L^{\perp}}(\hat\rho^{\eps}-S)P_{L^{\perp}},
\frac{1}{n}\sum_{j=1}^n \xi_j P_{L^{\perp}}X_jP_{L^{\perp}}\biggr\rangle
\biggr|
\leq \|P_{L^{\perp}}(\hat\rho^{\eps}-S)P_{L^{\perp}}\|_1
\biggl\|\frac{1}{n}\sum_{j=1}^n \xi_j P_{L^{\perp}}X_jP_{L^{\perp}}\biggr\|.
$$
By Proposition \ref{bernstein_A}, it is easy to show that 
with probability at least $1-e^{-t},$
$$  
\biggl\|\frac{1}{n}\sum_{j=1}^n \xi_j 
P_{L^{\perp}}X_j P_{L^{\perp}}\biggr\|\leq 
$$
$$
C\biggl[\sigma_{\xi} \|{\mathbb E}X^2\|^{1/2} 
\sqrt{\frac{t+\log (2m)}{n}}\bigvee
\|\xi\|_{\psi_2}\Bigl\|\|X\|\Bigr\|_{\psi_2} 
\log\biggl(\frac{\|\xi\|_{\psi_2}\Bigl\|\|X\|\Bigr\|_{\psi_2}}{\sigma_{\xi} \sigma_X}\biggr)\frac{t+\log (2m)}{n}\biggr].
$$
We replace $\sigma_X, \|{\mathbb E}X^2\|^{1/2}$ and 
$\Bigl\|\|X\|\Bigr\|_{\psi_2}$ by an upper 
bound $c\sqrt{m}$ (see Proposition \ref{subgauss}) which yields a simplified inequality
$$  
\biggl\|\frac{1}{n}\sum_{j=1}^n \xi_j 
P_{L^{\perp}}X_j P_{L^{\perp}}\biggr\|\leq 
C
\biggl[\sigma_{\xi} 
\sqrt{\frac{m(t+\log (2m))}{n}}\bigvee
\|\xi\|_{\psi_2}
\log\biggl(\frac{\|\xi\|_{\psi_2}}{\sigma_{\xi}}\biggr)\frac{\sqrt{m}(t+\log (2m))}{n}\biggr].
$$
Hence, with probability at least $1-e^{-t},$
$$
\biggl|\biggl\langle P_{L^{\perp}}(\hat\rho^{\eps}-S)P_{L^{\perp}},
\frac{1}{n}\sum_{j=1}^n \xi_j P_{L^{\perp}}X_j P_{L^{\perp}}\biggr\rangle
\biggr|\leq
$$
$$ 
C\|P_{L^{\perp}}(\hat\rho^{\eps}-S)P_{L^{\perp}}\|_1
\biggl[\sigma_{\xi} 
\sqrt{\frac{m(t+\log (2m))}{n}}\bigvee
\|\xi\|_{\psi_2}
\log\biggl(\frac{\|\xi\|_{\psi_2}}{\sigma_{\xi}}\biggr)\frac{\sqrt{m}(t+\log (2m))}{n}\biggr].
$$
The remaining term
$
\frac{1}{n}\sum_{j=1}^n \xi_j 
\langle \hat \rho^{\eps }-S,{\cal P}_L X_j\rangle 
$
is bounded exactly as in Step 3 of the proof of Theorem \ref{oracle} with the use of Adamczak's (2008)
version of Talagrand's concentration inequality.
This leads to the following bound: with probability at least $1-e^{-t},$
$$
\biggl|\frac{1}{n}\sum_{j=1}^n \xi_j 
\langle \hat \rho^{\eps }-S,{\cal P}_L X_j\rangle \biggr|
\leq 
C\biggl[\sigma_{\xi}\beta(L)
\|\hat \rho^{\eps}-S\|_{L_2(\Pi)} 
\sqrt{\frac{m r}{n}}+
\sigma_{\xi}  
\|\hat \rho^{\eps}-S\|_{L_2(\Pi)} 
\sqrt{\frac{\tau_n}{n}}+ \|\xi\|_{\psi_2}\frac{\sqrt{m} \tau_n\log n}{n}
\biggr],
$$
where $\tau_n=t+\log \log_2(2n).$

\qed


For simplicity, we state the next corollaries 
(similar to corollaries \ref{oracle-1} and \ref{oracle-Gibbs}) 
only in the case of 
subgaussian isotropic design. Recall that in this case $\|\cdot\|_{L_2(\Pi)}=\|\cdot \|_2$ and $\beta(L)=1.$

\begin{corollary} 
\label{oracle-2SG}
There exist numerical constants $C>0, c>0$ such that the following 
holds. For all $t>0$ and $\lambda>0$ 
such that $\tau_{n}\leq c\lambda^2 n,$ 
for all sufficiently large $D>0$ and for $\eps=D\eps_{n,m},$
for all matrices $S\in {\cal S}$ of rank $r,$
with probability at least $1-e^{-t},$ 
\begin{eqnarray}
\label{bound_6'''A-SG}
&&
\|\hat \rho^{\eps}-\rho\|_{L_2(\Pi)}^2
\leq (1+\lambda)\|S-\rho\|_{L_2(\Pi)}^2+
\frac{C}{\lambda}
\biggl[
D^2\biggl(\sigma_{\xi}^2 \frac{r m t_{m}}{n}\bigvee
c_{\xi}^2\frac{r m t_{m}^2}{n^2}\biggr)\log^2 (m n)
\bigvee 
\nonumber
\\
&&
\sigma_{\xi}^2 \frac{\tau_n}{n} 
\bigvee
(c_{\xi} \vee \sqrt{m}) \frac{\sqrt{m} t_{n,m}}{n}\biggr].
\end{eqnarray}
\end{corollary}

\begin{corollary} 
\label{oracle-Gibbs-3}
There exists numerical constants $C>0,c>0$ such that the following 
holds. For all $t>0$ and for all $\lambda>0$  
such that $\tau_n\leq c\lambda^2 n,$ 
for all sufficiently 
large $D$ and for $\eps=D\eps_{n,m},$
for all Hermitian matrices $H$ and for all $r\leq m,$  
with probability at least $1-e^{-t},$ 
\begin{eqnarray}
\label{bound_6-Gibbs-3}
&&
\|\hat \rho^{\eps}-\rho\|_{L_2(\Pi)}^2
\leq (1+\lambda) \|\rho_{H}-\rho\|_{L_2(\Pi)}^2+
\frac{C}{\lambda}
\biggl[\delta_r^2(H)\bigvee 
D^2\biggl(\sigma_{\xi}^2 \frac{\Gamma_r m t_{m}}{n}\bigvee
\nonumber
\\
&&
c_{\xi}^2\frac{\Gamma_r m t_{m}^2}{n^2}\biggr)
\bigvee 
\sigma_{\xi}^2 \frac{m r+\tau_n}{n} \bigvee 
(c_{\xi} \vee \sqrt{m}) \frac{\sqrt{m} t_{n,m}}{n}
\biggr].
\end{eqnarray}
\end{corollary}

In a special case of Gaussian noise, the bounds of the above corollaries 
can be simplified since in this case 
$c_{\xi}\leq c\sigma_{\xi}$ for some numerical 
constant $c.$ In particular, Corollary \ref{oracle-2SG} immediately implies the bound of Theorem \ref{intro_th_2} in the Introduction.
Both bounds of Theorem \ref{intro_th_0} follow from 
theorems \ref{simple-SG} and \ref{oracle-SG}.

}

\end{document}